\documentclass{article}

\AtBeginDocument{%
	\providecommand\BibTeX{{%
			\normalfont B\kern-0.5em{\scshape i\kern-0.25em b}\kern-0.8em\TeX}}}

\usepackage{arxiv}
\usepackage{graphicx}
\graphicspath{{./figs/}}
\usepackage{color}
\usepackage{array}
\usepackage{multirow}
\usepackage{threeparttable}
\usepackage{stfloats}
\usepackage{amssymb, amsmath}
\usepackage{psfrag}
\usepackage{csquotes}
\usepackage{fancyhdr} % for header and footer
\usepackage[symbol]{footmisc}

\usepackage{IEEEtrantools}
\newcommand{\phiplus}{$\phi^+(x)$}
\newcommand{\phiminus}{$\phi^-(x)$}

\newcommand{\Phipm}{$\Phi^{\pm}(x)$}
\newcommand{\LULP}{$\text{ULP}_{\text{LNS}}$}
\newcommand{\ULP}{$\text{ULP}$}
\newcolumntype{C}[1]{>{\centering\let\newline\\\arraybackslash\hspace{0pt}}m{#1}}

\title{Low Precision Logarithmic Number Systems: Beyond Base-2}

\author{Syed Asad Alam, James Garland and David Gregg%\thanks{Lero, Trinity College Dublin, College 
 \\
 School of Computer Science and Statistics\\
 Lero, Trinity College Dublin\\
 College Green, Dublin 2\\
 Dublin, Ireland\\
 \texttt{\{syed.asad.alam,jgarland,dgregg\}@tcd.ie}}

\begin{document}
 
 \maketitle	

\begin{abstract}
Logarithmic number systems (LNS) are used to represent real numbers in many applications using a 
constant base raised to a fixed-point exponent making its distribution exponential. This greatly 
simplifies hardware multiply, divide and square root. LNS with base-$2$ is most common, but in this 
paper we show that for low-precision LNS the choice of base has a significant impact.

We make four main contributions. First, LNS is not closed under addition and subtraction, so the 
result is approximate. We show that choosing a suitable base can manipulate the distribution to 
reduce the average error. Second, we show that low-precision LNS addition and subtraction can be 
implemented efficiently in logic rather than commonly used ROM lookup tables, the complexity of 
which can be reduced by an appropriate choice of base. A similar effect is shown where the result 
of arithmetic has greater precision than the input. Third, where input data from external sources 
is not expected to be in LNS, we can reduce the conversion error by selecting a LNS base to match 
the expected distribution of the input. Thus, there is no one base which gives the global optimum, 
and base selection is a trade-off between different factors. Fourth, we show that circuits realized 
in LNS require lower area and power consumption for short word lengths. 
 
\end{abstract}

%\begin{CCSXML}
%	<ccs2012>
%	<concept>
%	<concept_id>10010520.10010553.10010562</concept_id>
%	<concept_desc>Computer systems organization~Embedded systems</concept_desc>
%	<concept_significance>500</concept_significance>
%	</concept>
%	<concept>
%	<concept_id>10002950.10003714.10003715</concept_id>
%	<concept_desc>Mathematics of computing</concept_desc>
%	<concept_significance>300</concept_significance>
%	</concept>
%	</ccs2012>
%\end{CCSXML}

%\ccsdesc[500]{Computer systems organization~Embedded systems}
%\ccsdesc{Hardware~Integrated circuits~Logic circuits~Arithmetic and datapath circuits}
% I commented the next line cos it causes latex to fail
%\ccsdesc[300]{Mathematics of computing}


\maketitle

\section{Introduction}
\label{sec:intro}

%%% A gentle introduction focusing on technologies
As the world moves towards autonomous systems in every field, there is a need for small embedded 
devices that can interact with the environment on its edge rather than offloading computations to 
the cloud. These devices often operate on physical-world data such as audio, video and input from 
various sensors, using computationally intensive techniques, like filtering, facial recognition, 
image segmentation and object detection. Such techniques are based on signal processing or deep 
neural networks (DNN). The ability to execute these applications on small resource-constrained 
embedded devices is a major milestones in truly enabling technologies like driverless cars, smart 
cities and smart mobile devices \cite{Satyanarayanan2017J}. Given the increasing quantities of data 
processed by smart devices and autonomous systems, computing at the edge is essential. Sending data 
back to the cloud for processing requires large amounts of energy, and may become a bottleneck that 
fails to meet latency requirements for real-time embedded and autonomous systems while also 
stretching the limited processing, memory and/or battery capacity of embedded devices 
\cite{Shi2016_J, Chen_2019J}.

%%% Key aspects for enabling DNNs on embedded devices --> Making it generic and not specific to DNNs
%Processing large amounts of data at the edge is a challenge for embedded devices that may
%have limited processing, memory, and/or battery capacity for applications listed 
%earlier. 
One possible solution is a logarithmic number system (LNS) which represents real numbers as a 
fixed-point exponent. Compared to floating point, LNS dramatically simplifies the hardware needed 
for multiplication, division, and square root. Multiplication and division become fixed-point 
addition and subtraction of the exponents, which can be implemented with a simple integer adder 
\cite{asad_lns2014,Kouretas_2013J}. Square root is computed by dividing the exponent by two, which 
can be implemented by simply dropping the least-significant bit of the fixed-point exponent. Among 
these basic operations, multiplication occur very frequently in systems like neural networks and 
finite-impulse response (FIR) filters and reduction of its complexity will have a major impact on 
the overall reduction of area and power. The downside of LNS is that addition and subtraction are 
non-trivial. A common strategy to implement LNS adders/subtracters use ROM look-up 
tables, which in the simple case are exponential in the size of LNS word. There is an extensive 
literature on compacting these tables for large word sizes \cite{kenny_2008, Coleman2008J, 
Mansour2015J}.

Smaller word sizes can reduce the memory footprint of data, and the complexity of 
arithmetic for a variety of number systems \cite{Feng2010, Rastegari2016X, Capotondi2020J, 
Garofalo2019J, Bruschi2020C, ChungBrainwaveMicrosoft_2018M} which 
will have a significant impact on the ability to deploy systems in resource constrained embedded 
devices deployed on the edge of networks. Further fixed point number systems with very short word 
lengths have been proposed in the literature for a variety of signal processing applications 
\cite{Lawrence1980,asad_lns2014} while shorter floating and fixed point numbers have also been used 
for neural networks \cite{Hubara2017J,Wang2018trainingC,Wang_2018J}.%,sun2019hybridC}. 

%Similarly, short word 
%length LNSs have also been proposed for various other applications 
%\cite{Lee2017C,Vogel2018C,Miyashita2016X,asad_lns2014}.

An important aspect of any number system is the distribution of the representable numbers. LNS have 
a distribution that is non-uniform. Each number is a constant multiple of the previous one,
resulting in an exponential distribution. Fully one half of LNS numbers fall between $-1$ and $1$, 
and representable numbers become increasingly sparse thereafter. In applications like FIR filters, 
the coefficient space lies within the range of $\pm1$ \cite{asad_lns2014} while the distribution of 
weights and biases in DNNs is also typically clustered within a small range \cite{Lee2017C}. Thus, 
LNS is suitable for use in a number of different applications \cite{Miyashita2016X, asad_lns2014, 
Kouretas_2013J}.

The exponential distribution of LNS depends on two factors: the base and the number of 
integer/fractional bits in the fixed-point exponent. In existing literature, the base is almost 
always a power of two, typically $2$, $\sqrt{2}$, or $\sqrt[4]{2}$. However, there is no
systematic study of different bases in the literature, perhaps because within sensible ranges the 
choice of base has little impact on LNS with large word size. However, we show that for 
low-precision LNS the base has a major impact. We make the following contributions.

\begin{itemize}
	\item We demonstrate that LNS with different bases have different arithmetic rounding errors which 
	is significant for low-precision LNS and can be reduced by selecting an appropriate base. 
	\item We optimize hardware for low-precision LNS arithmetic by selecting a base with favourable 
	add and subtract truth tables. We achieve large savings in circuit area and depth for LNS adder 
	and subtractor circuits by implementing the tables using logic gates as compared to ROM based 
	implementation.%, and for hybrid 	logic/lookup-table implementations.
 \item We show that where input data must be converted from another number format to LNS, base 
 selection can enable a more accurate approximation of the original data distribution, where 
 accuracy is measured by the representation error.%. We measure accuracy by measuring the 
 %representation error.
 \item We also define arithmetic and representation error in a multiplicative sense in the real 
 domain which allows us to describe the error in bits.
	\item We analyze mixed precision tables where the values stored in addition and subtraction tables 
	have more precision than their inputs and show that realizing these tables results in a reduced 
	error without an exponential increase in the hardware cost to realize them.
 \item We synthesis adders, subtractors and FIR filters for different word lengths using 
 different number systems and show that circuits LNS based circuits 
 achieve lower area and power consumption as compared to its counterparts for short word lengths.
\end{itemize}

The paper is organized as follows: Section~\ref{sec:related_work}
describes related work. Section~\ref{sec:lns} gives a brief
introduction of logarithmic number system, the effect of base aliasing
and of rounding of numbers in LNS. Section~\ref{sec:quant_float}
discusses the quantization of floating point numbers into LNS.
Section~\ref{sec:lns_addsub} describes LNS addition and subtraction
using tables and approximation errors in these two operations. The
synthesis of addition and subtraction tables is described and
evaluated in Section~\ref{sec:lns_addsub_synthesis} while
Section~\ref{sec:tradeoff_lnsaddsub} presents a discussion on the
trade-offs in selecting the best base based on various parameters.
Mixed-precision tables are the focus in
Section~\ref{sec:high_precision_outputs} while
Section~\ref{sec:lns_v_fixedfloat} discusses the hardware cost of
realizing LNS adder/subtractor circuits and compared to fixed point
multipliers and floating point adders and multipliers. The discussion
of hardware cost of real word circuits is extended in
Section~\ref{sec:fir_filter} for FIR filters using the different
number systems and word length.
%Section~\ref{sec:conc} concludes the paper.

\section{Related Work}
\label{sec:related_work}

%With these properties,
LNS has been evaluated as an alternative to fixed and floating point numbers in a number of 
applications. It provides better dynamic range and round-off noise performance as compared to 
fixed-point \cite{king_1971} and floating point \cite{arnold_2002} in some signal processing 
applications. Paliouras and Stouratis \cite{pal_2000} show that LNS can lead to lower switching 
activity leading to reduced power consumption while also reducing representational error relative 
to fixed point. Basetas et al. in \cite{baset_2007} show that LNS requires a reduced word length 
when implementing finite impulse response (FIR) and infinite impulse response (IIR) filters. Alam 
and Gustafsson have also shown \cite{asad_lns2014} that LNS improves the 
approximation error by around $20\%$ when designing FIR filters as compared to a 
fixed point design. Coleman proposed a $32$-bit logarithmic processor and concludes that LNS 
provides faster execution, more accuracy and reduced architectural complexity as compared to single 
precision $32$-bit floating point, provided that the adder/subtractor has a latency less than that 
of a $32\times32$-bit multiplier. Chandra \cite{satish_1998} shows an order of magnitude 
improvement in error-to-signal variances for round-off noise in FIR filters using LNS over floating 
point number representation.

%studies the effect of 
%round-off noise in FIR filters using LNS and compares it against that of floating point number 
%representation and shows an order of magnitude improvement in error-to-signal variances.

LNS has also been explored for deep neural networks, especially for very short 
word lengths. The LogNet inference engine, proposed by Lee et al. \cite{Lee2017C}, uses $4$- and 
$5$-bit LNS encoding to represent weights and achieves improved results as compared to fixed point 
encoding without retraining while keeping the activations in floating point. They also encode the 
activations in LNS with $3$- and $4$-bits while keeping the weights in floating point and show 
negligible to zero degradation in performance as compared to floating-point. For 
e.g., the top-5 accuracies (without training) using logarithmic encoding in base-$2$ with $4$ and 
$5$ bits is $73.4\%$ and $74.6\%$, respectively, for AlexNet and $85.2\%$ and $86.0\%$, 
respectively, for VGG16. The corresponding accuracy using $32$-bit floating point numbers is 
$78.3\%$ and $89.8\%$ for the two nets, respectively and much better than corresponding fixed-point 
performance. This number further improves to a maximum of $85.2\$$ upon training and logarithmic 
encoding also has more sparse weights ($35.8\%$ against $21.2\%$) allowing for reduced memory 
requirements. However, although not clearly stated, the exponent for LNS is represented using 
integers only and the impact of representing the exponent using fixed-point numbers is not explored.

Miyashita et al. also use an integer representation for the exponent while 
quantizing using $3$- and $4$-bit encodings for activations and $5$-bit for weights 
\cite{Miyashita2016X}. For AlexNet, Miyashita et al. show that a $4$-bit encoding achieves only a 
$1.4\%$ drop in the top-5 accuracy for AlexNet and no drop for VGG16. The corresponding drop in 
accuracy for $3$-bit logarithmic encoding is also negligible. Miyashita et al. also report a 
reduction in the memory requirements for AlexNet and VGG16 using LNS encoding of $85.7\%$ and 
$77.9\%$, respectively. They have also considered bases of $\sqrt{2}$ and $\sqrt[4]{2}$ for 
representation and find that the total quantization error on the weights for $ \sqrt{2} $ is about 
$2\times$ smaller than base-$2$ LNS. However, as we show in Section \ref{subsec:lns_basealias}, 
bases $2$, $\sqrt{2}$ and $\sqrt[4]{2}$ are simple aliases of one another, where the binary point 
of the fixed point exponent is moved left. Vogel et al.\cite{Vogel2018C}, present a low word length 
based quantization method for LNS to be used in a neural network and have shown LNS to achieve 
$22.3\%$ lower power as compared to an $8$-bit fixed point based design. Arnold et al. 
\cite{Arnold2020C} present an approach to implement back propagation using tableless LNS ALU based 
on modified Mitchell's method \cite{Mitchell1962J} and achieve one third reduction in the hardware 
resources as compared to a conventional fixed-point implementation.

However, to the best of authors' knowledge, a comprehensive and systematic study 
of logarithmic number systems (LNSs) with respect to different bases and error measurement is 
missing in the literature. Also missing is the study of the impact of choice of base on short word 
length systems and providing a study of the trade-off involved in selecting bases as a function of 
the word length while also analyzing the impact of LNS on systems like FIR filters which are 
heavily dependent on the operations of multiplications and additions. These important aspects form 
the bulk of the contribution of this work.

\section{Logarithmic Number System (LNS)}
\label{sec:lns}

A number in the logarithmic number system is represented using the zero-bit, sign-bit and the 
actual logarithmic value, represented as a triplet ($\mathcal{X}$), as follows \cite{earl_1975}:

\begin{equation}
\mathcal{X} = (z_x,s_x,m_x),
\label{eq:log}
\end{equation}
where $z_x$ indicates whether $\mathcal{X}$ is zero, $s_x = sign(X)$ and $m_x = \log_b|X|$. The 
exponent (or logarithmic) value $m_x$ is typically represented using two's complement fixed-point 
representation, with $i$ and $f$ integer and fractional bits (mathematically represented as 
$Q(i,f)$) \cite{asad_lns2014, Coleman2008J}. The base-$b$ of an LNS number influences the 
representation capabilities and complexity of computations and conversions \cite{stouraitis_02}. 
The real value zero cannot be expressed in the form $b^m, b \neq 0$, and therefore LNS commonly
uses an additional bit to represent zero. Some LNS instead use the most negative possible $m_x$ as 
a special value to represent zero \cite{kenny_2008}.

\subsection{Bases and Base Aliasing in LNS}
\label{subsec:lns_basealias}

The base $b$ is a key parameter of the LNS. For a base $b > 1$, each LNS value is exponentially 
larger than the previous one. However, the ratio between successive LNS values is determined not 
just by the base, but also by the number of fractional bits in the fixed point format $Q(i,f)$. The 
least significant bit of the fixed-point exponent has the value $2^{-f}$ in the \textit{log domain}
of the LNS number which we define as the \textit{unit of least-precision} in the log domain 
(\LULP).  Increasing the log domain value by one \LULP~corresponds to multiplying the LNS number by 
$b^{2^{-f}}$ in the \textit{real domain} which we simply define as \ULP. We further define an 
integer radix form, $r$ of the LNS base, where $r = b^{2{^-f}}$. The ratio between successive 
values of the LNS  in the real domain is $r$.

Note that the LNS with $Q(i, f)$ and base $b$ contains the same set of values as the LNS with 
$Q(i+f, 0)$ and base $r$. Thus, the effect of changing the base of an LNS may be equivalent to 
simply moving the binary point in the fixed-point representation $m_x$.  For example, Miyashita et 
al. \cite{Miyashita2016X} explored the effect of varying the LNS base by experimentally evaluating 
the bases $2$, $\sqrt{2}$ and $\sqrt[4]{2}$. However:

\begin{IEEEeqnarray}{rCl}
	\sqrt{2}^{i_2i_1i_0.f_1f_2f_3f_4} & = & 2^{\frac{{i_2i_1i_0.f_1f_2f_3f_4}}{2}}\\ \nonumber
									  & = & 2^{{i_1i_0.f_1f_2f_3f_4f_5}}
	\label{eqn:base_alias}
\end{IEEEeqnarray}

In other words, LNS with base-$\sqrt{2}$ and exponent of the form $Q(i,f)$ contains exactly the 
same set of values as base-$2$ LNS with $Q(i-1, f+1)$, and base-$\sqrt[4]{2}$ with $Q(i+1, f-1)$. 
By selecting the position of the fixed point, standard base-$2$ LNS can represent implicit bases 
such as $\sqrt[8]{2}$, $\sqrt[4]{2}$, $\sqrt[2]{2}$, $ 2 $, $ 4 $, $ 8 $, etc. By exploring only 
bases that can 
be already represented using base-$2$, Miyashita et al. \cite{Miyashita2016X} investigated only the 
effect of moving the binary point in $m_x$, not in exploring arbitrary bases. Remarkably, despite a 
very extensive literature on LNS over several decades, we have been unable to find a statement
of this simple identity in the literature.

Thus, when choosing a base for the LNS, the important factors are the total number of bits in 
$m_x$, that is $i+f$, and the integer radix form of the base, $r$. In principle, fixed-point
exponents are entirely unnecessary, and we can simply use an integer exponent and a suitable base 
$b = r$. However, practical values of $r$ are typically not very convenient numbers to refer
to. Therefore, we maintain the existing tradition of representing $m_x$ as a fixed-point number, 
rather than an integer. We explore bases in the range $\sqrt{2}$ to $2$ which covers a range of 
values separated by one binary point position in the format of $m_x$.

\subsection{Rounding LNS numbers}
\label{subsec:lns_rounding}

LNS values have a finite number of digits, and therefore cannot represent all real numbers exactly. 
Where a real number arising from arithmetic cannot be represented exactly in the LNS, we normally 
round to a nearby representable value. We focus on rounding to the nearest representable LNS 
number. Rounding to nearest is a common strategy in floating point number systems; in 
non-exceptional cases the FP rounding error is bounded by half of one FP mantissa unit of least 
precision. 

A problem with rounding to nearest in LNS is that there are two different and incompatible possible 
approaches. One option is to round in the LNS domain, i.e., to round the LNS exponent to the nearest
representable fixed-point exponent $m_x$. Another option is to round in the real domain, that is to 
consider the LNS number and its nearest neighbours to the real domain, and choose the rounding that 
minimizes the real domain error.

The two approaches give slightly different answers. Consider the case of rounding the fixed-point 
\textit{binary} exponent value $0.1001$ in a base-$2$ LNS with only integer digits. The choice is 
between rounding this value downwards to an exponent of 0, or upwards to an exponent of 1. If we 
round in the log domain, it is clear that $0.1001$ is greater than the mid-point of $0.1$ and we 
should therefore round upwards to 1.  On the other hand, if we convert these numbers to the real 
domain, then we are choosing between rounding downwards to $2^0 = 1$ or upwards to $2^1 = 2$. The 
value of binary fixed point exponent $0.1001$ is $2^{\frac{9}{16}} \approx 1.44043$, which is 
clearly closer to 1 than 2. Thus, if we round in the real domain, we should round this number 
downwards.

In addition to rounding in the real or log domain, one needs to decide whether to measure absolute 
or relative error.  For fixed point numbers, one either uses the absolute error between the 
original value and the rounded one or the relative error, as given here:

\begin{equation}
 e_a = \hat{m_x} - m_x
 \label{eq:abserr}
\end{equation}

\begin{equation}
 e_r = \frac{e_a}{m_x}
 \label{eq:relerr}
\end{equation}
where $e_a$ and $e_r$ represent the absolute and relative error, respectively.

The measure of error used in literature has typically been the relative error but in the real 
domain \cite{Chandra1998J}, as given by:

\begin{equation}
 e_{rr} = \frac{b^{\hat{m_x}} - b^{m_x}}{b^{m_x}}
\label{eq:add_rel_err}
\end{equation}
where $e_{rr}$ is the relative error in the real domain and $b$ is the base. However, recall that 
in LNS adding one \LULP~to a number in the log domain corresponds to multiplying the real domain 
value by the integer radix form, $r$, of the base, given by $b^{2^{-f}}$ and defined as \ULP~in 
Section~\ref{subsec:lns_basealias}. Thus, rounding in the log domain minimizes the 
\textit{multiplicative error} of the real domain value, that is the factor multiple between the 
exact real-domain value and the rounded real domain value. On the other hand, rounding in the real
domain minimizes the real domain \textit{additive relative error}, as given in Equation 
(\ref{eq:add_rel_err}), that is the amount that must be added to or subtracted from the exact real 
domain value to get the rounded real domain value. We define the multiplicative error in the real 
domain as: 

%This is the reason that ULP in the real domain for quantization in LNS is defined as $b^{2^{-f}}$. 

\begin{equation}
 e_m = \frac{b^{\hat{m_x}}}{b^{m_x}} = b^{\hat{m_x} - m_x}
 \label{eq:mut_error}
\end{equation}
where $b$ is the base, as usual. The range of this error, ignoring round to zeros, is given by: 

\begin{equation}
 \frac{1}{\text{ULP}_h} \le e_m \le \text{ULP}_h
 \label{eq:mult_error_range}
\end{equation}
where $\text{ULP} = b^{2^{-f}}$ and $\text{ULP}_h = b^{2^{-f-1}}$ with $b$ being the base and $f$ 
being the number of fractional bits. Since this error measure is based on a ratio, it makes more 
sense to calculate the geometric mean of it, as given by:

\begin{equation}
 \mu_g = \left(\prod_{i=1}^{n}e_{m,i}\right)^{\frac{1}{n}}
 \label{eq:geo_mean}
\end{equation}
where $e_{m,i}$ is the individual multiplicative relative error in the real domain given by 
Equation~\ref{eq:mut_error} and $n$ is the total number of values over which the mean is 
calculated. 

The geometric mean corresponds to the exponential of the arithmetic mean of logarithms given by:

\begin{IEEEeqnarray}{rCl}
 \left(\prod_{i=1}^{n}e_{m,i}\right)^{\frac{1}{n}} & = & 
 b^{\left[\frac{1}{n}\sum_{i=i}^{n}log_b(e_{m,i})\right]} \\\nonumber
 & = & b^{\left[\frac{1}{n}\sum_{i=i}^{n}(log_b(b^{\hat{m_{x,i}} - m_{x,i}}))\right]} \\\nonumber
 & = & b^{\left[\frac{1}{n}\sum_{i=i}^{n}\hat{m_{x,i}} - m_{x,i}\right]}
 \label{eq:gm_am}
\end{IEEEeqnarray}

This implies that calculating the geometric mean in the real domain corresponds to calculating 
arithmetic mean of the absolute error of the exponents which corresponds to rounding in the log 
domain.

Furthermore, rounding in the log domain guarantees that in non-exceptional cases the absolute 
rounding error is bound by half of one \LULP. However, a problem arises when rounding to 
real-domain zero, which has a log domain value roughly equal to $-\infty$. Thus, when rounding a 
non-zero value to zero, the log-domain error is infinite. If our rounding goal is to minimize the 
log domain error, we would never round a value to real domain zero. However, never rounding to zero 
does not make sense, at least to us. Therefore, we round values of less than half of one log domain 
\LULP~down to zero.

%Nonetheless, we round all values less than $\sqrt{r}$ down to zero.

Rounding in the real domain eliminates the rounding-to-zero anomaly of an infinite error. But 
rounding to zero causes the real-domain relative errors to be $100\%$. In all other cases of 
rounding, the real domain relative error is bound by $\sqrt{r}-1$ for all $r > 1$. Despite a very 
extensive literature on LNS over decades, we have not been able to locate a similar discussion of 
log-domain versus real-domain rounding. In the remainder of the paper we round in the log domain, 
which allows us to present errors in $\textit{ULP}_{\textit{LNS}}$ (which equals $2^{-f}$ in the 
log domain), but rounding in the real domain gives similar results. Furthermore, when presenting 
results on the average absolute error in the log domain, we report values that underflow to zero 
--- and thus result in an infinite log-domain error --- separately from from values that
do not underflow.

\section{Conversion of non-LNS inputs}
\label{sec:quant_float}

When computing with LNS, it is common that the original inputs to the program are in a different 
format, such as floating point (FP). We normally convert inputs to LNS once, and perform all 
subsequent computation in LNS. Any LNS and FP format contain different sets of real values. Thus, 
when we convert input values to LNS, we normally round each input to the nearest LNS number. Note 
that the rounding error for converting inputs tend to be cumulatively smaller than arithmetic 
rounding errors, because most algorithms perform multiple arithmetic operations for each input. In 
this section we show that we can reduce the conversion rounding error, as compared to base-$2$ LNS,
by selecting another base that better captures the distribution of the input data.

We consider a simple case of input data conversion, where the input consists of a short length FP 
input type which we convert to a logarithmic number format with the same number of bits. Each FP 
number has a sign bit, $i$ exponent bits, and $f$ mantissa bits, and we convert to an LNS type with 
a sign bit, and an exponent with $i$ integer bits and $f$ fractional bits. Very low precision types 
have been used for both FP and LNS, even as low as $3$-bits \cite{Miyashita2016X,Lee2017C}. In 
light of these contributions it is important to analyze the effect of quantization for shorter word
lengths. The word length formats used for our work are given in Table~\ref{tab:number_format}.
 
\begin{table}
	\centering
	\caption{Number formats used for representing LNS numbers.}
	\label{tab:number_format}
	\begin{threeparttable}
  \begin{tabular}{C{1cm}C{1cm}C{0.75cm}c|C{1cm}C{1cm}C{0.75cm}c}
		\hline
		Word length\tnote{$\dagger$} & Int. bits & Frac. bits & Notation&
  Word length\tnote{$\dagger$} & Int. bits & Frac. bits & Notation\\
		\hline
		$5$  & $2$ & $2 $ & Q(2,2)  & $6$   & $2$ & $3 $ & Q(2,3) \\
		$7$  & $3$ & $3 $ & Q(3,3)  & $8$   & $4$ & $3 $ & Q(4,3) \\
  $9$  & $4$ & $4 $ & Q(4,4)  & $10$  & $4$ & $5 $ & Q(4,5) \\
  $11$ & $4$ & $6 $ & Q(4,6)  & $12$  & $4$ & $7 $ & Q(4,7) \\
  $13$ & $4$ & $8 $ & Q(4,8)  & $14$  & $4$ & $9 $ & Q(4,9) \\
  $15$ & $4$ & $10$ & Q(4,10) & $16$ & $5$ & $10$ & Q(5,10)\\
  \hline
	\end{tabular}
 \begin{tablenotes}
  \item[$\dagger$] $1$-bit reserved for sign.
 \end{tablenotes}
 \end{threeparttable}
\end{table}

In our experiments we assume that all FP numbers of the appropriate input type are equally likely 
to appear. In addition to considering base-$2$ LNS, we consider around $ 587 $ evenly distributed 
bases between $\sqrt{2}$ and $2$. When computing the rounding error in conversion, we use the 
arithmetic mean of the absolute error of the exponents (in the log domain) as discussed in 
Section~\ref{subsec:lns_rounding} and given in Equations~\ref{eq:abserr} and \ref{eq:gm_am}. 

\begin{figure}
	\centering
	\includegraphics[scale=0.6]{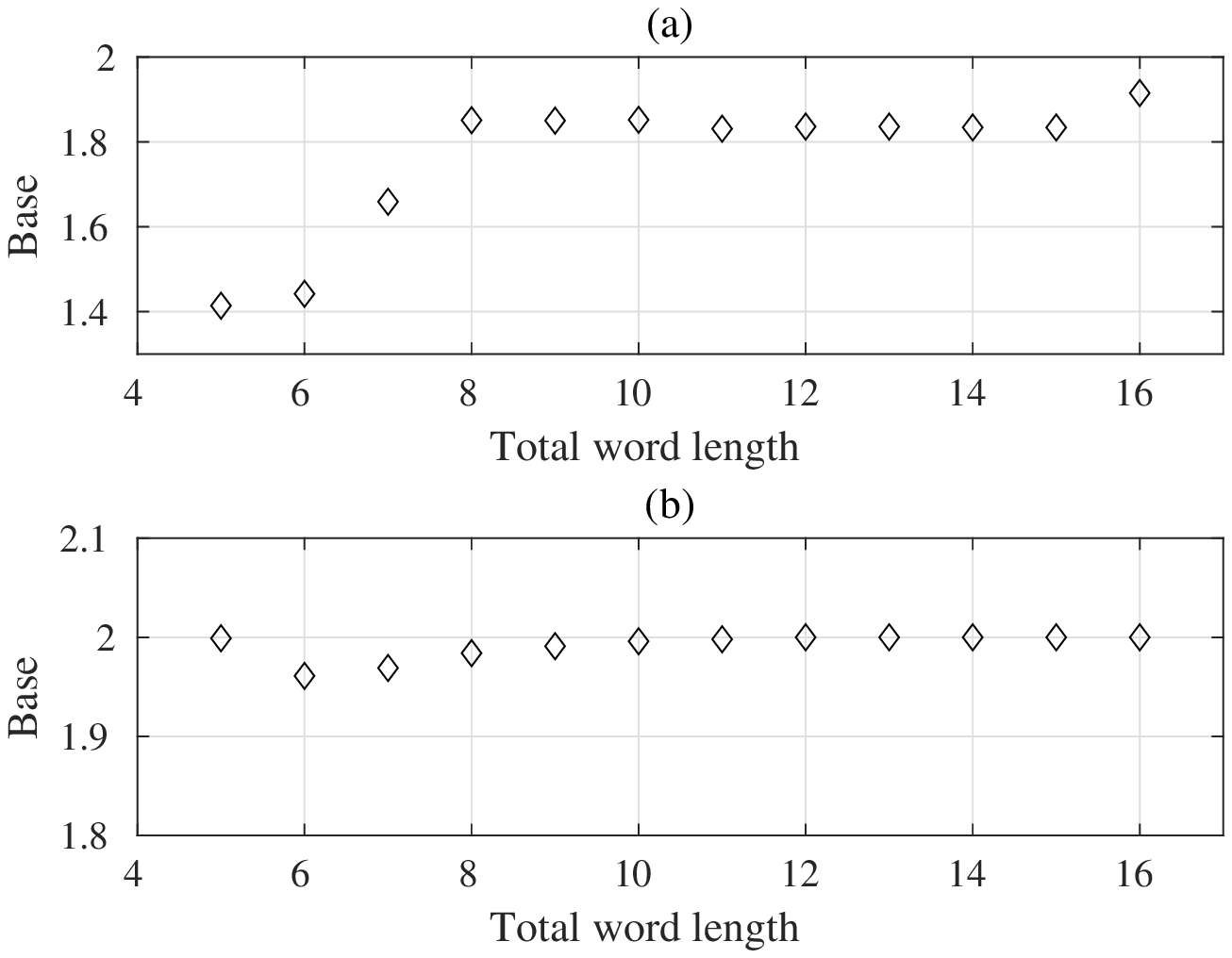}
	\caption{Bases with minimum quantization error for various word lengths for (a) scaled and (b) unscaled LNS data.}
	\label{fig:lnsquant_scaled_mindata}
\end{figure}

The problem of rounding a non-zero FP number to zero during conversion is the same as stated 
earlier in Section~\ref{subsec:lns_rounding} and we exclude values that underflow or overflow
during rounding, as also stated earlier in Section~\ref{subsec:lns_rounding}. 
Fig.~\ref{fig:lnsquant_scaled_mindata}(b) shows that for the distribution of FP numbers, the LNS 
base that gives the lowest quantization error is typically close to $2$ for all word lengths given 
in Table~\ref{tab:number_format}.

Note, however, that each LNS base gives both a different \textit{range} of real-domain
values, and a different \textit{distribution} of those values, with bases close to 
base-$2$ having a range closest to a floating point distribution for a given word length. This is 
reflected in Fig.~\ref{fig:lnsquant_scaled_mindata}(b) where bases close to base-$2$ provide the 
least quantization error across all word lengths. We consider a second experiment where we apply a 
constant scaling factor to all LNS values, so that the range of LNS values is the same for all 
bases. With this scaling, different bases lead to different distributions of LNS values,
but the range is constant. We use a scaling factor that is the ratio of the maximum
floating point and LNS number for each base. Fig.~\ref{fig:lnsquant_scaled_mindata}(a)
shows the effect of keeping the range of LNS values constant through scaling with 
different bases resulting in distributions that can reduce the quantization error, compared to 
base-$2$. For very short word lengths of $5$ -- $7$, the optimal base for quantization error is 
much lower than base-$2$.

%%% Referring to Table 2
The average absolute representation error, as a percentage of \LULP~in the log domain and for a 
scaled LNS distribution, for a few select bases is shown in Table~\ref{tab:perc_phi_bases} along 
with the difference of other bases to the best base that gives the minimum representation error, in 
columns two and three. Other columns describe another type of error which will be explained shortly.

It can be seen from Table~\ref{tab:perc_phi_bases} that bases other
than base-$2$ can give significant reduction in error, with a minimum
$4\times$ improvement for $Q(2,2)$ over base-$2$ and a maximum of
close to three order of magnitude for larger word lengths. However,
the distribution and range of inputs may be domain dependent. A base
that results in a distribution of values that is suitable for one
domain might be a poor match for another domain. Thus, matching the
LNS base to an input distribution may be useful when designing a
hardware accelerator that targets a specific domain. But for
general-purpose LNS processing, it is probably not worthwhile to build
a base that is specific to one domain into the hardware LNS units. Our
experiments show that we can adapt the LNS base to reduce the rounding
when converting \textit{some} input distribution to LNS. However, any
reduction in conversion error may also be domain-dependent. In
contrast, our approach in the next section improves the accuracy of
arithmetic, and is less domain-dependent.

\begin{table}
 \footnotesize
 \centering    
 \begin{threeparttable}
  \caption{Percentage average relative error for various bases.}
  \label{tab:perc_phi_bases}
  \begin{tabular}{c|c|c|cc|cc}
   \hline
   \multirow{2}{*}{Number format} & \multirow{2}{*}{Base} & \multirow{2}{*}{Avg. abs. repr. error 
    (\%)\tnote{$\dagger$}} & \multicolumn{4}{c}{Avg. abs. arith. error (\%)} \\
   \cline{4-7} 
   & & & \phiplus & $N_z$\tnote{$\bigstar$} & \phiminus & $N_z$\tnote{$\bigstar$} \\
   \hline
   \multirow{4}{*}{Q(2,2)} 
   &$1.414$&$\mathbf{2.6888}$&$22.9985$	        &$0$&$20.5680$         &$4$ \\
   &$1.741$&$5.4437$         &$\mathbf{22.1167}$&$0$&$21.5063$         &$1$ \\
   &$1.417$&$2.8546$         &$22.7008$	        &$0$&$\mathbf{17.1544}$&$4$ \\
   &$2.0$  &$9.7655$	        &$26.0090$	        &$0$&$25.4259$         &$1$ \\
   \hline
   
   \multirow{4}{*}{Q(2,3)}          
   & $1.442$&$\mathbf{1.2682}$&$24.2746$         &$0$&$26.2439$          & $7$ \\
   & $1.999$&$21.6569$        &$\mathbf{22.9960}$&$0$&$21.1268$          & $7$ \\
   & $1.415$&$1.6456 $        &$24.8867$         &$0$&$\mathbf{19.9971}$ & $1$ \\
   & $2.0$  &$23.6369$        &$23.6633$         &$0$&$21.4604$          & $1$ \\
   \hline
   
   \multirow{4}{*}{Q(3,3)}
   & $1.659$& $\mathbf{1.6902}$&$25.4686$         & $0$&$24.7396$         &$1$ \\
   & $1.496$& $9.6850$         &$\mathbf{23.4968}$& $0$&$23.1648$         &$1$ \\
   & $1.802$& $3.5171$         &$25.0584$         & $0$&$\mathbf{22.1066}$&$1$ \\
   & $2.0$  & $23.6369$        &$25.3094$         & $0$&$23.7062$         &$1$ \\
   \hline
   
   \multirow{4}{*}{Q(4,3)}          
   & $1.851$&$\mathbf{1.6070}$&$25.2744$         &$0$&$24.6680$         &$1$ \\ 
   & $1.730$&$5.0651	        $&$\mathbf{24.8420}$&$0$&$25.1975$          &$1$ \\ 
   & $1.802$&$2.4998	        $&$25.0584$         &$0$&$\mathbf{22.1066}$&$1$ \\ 
   & $2.0$  &$23.6369        $&$25.3094$         &$0$&$23.7062$         &$1$ \\ 
   \hline
   
   \multirow{4}{*}{Q(4,4)}                                       
   & $1.850$&$\mathbf{1.2107}$&$25.2267$         &$0$&$25.6736$        &$1$ \\
   & $1.697$&$6.1957	$        &$\mathbf{25.0533}$&$0$&$24.0765$         &$1$ \\
   & $1.973$&$6.5001	$        &$25.2602$         &$0$&$\mathbf{23.4122}$&$1$ \\
   & $2.0$  &$14.9733$        &$25.5656$         &$0$&$24.7483$         &$1$ \\
   \hline
   
   \multirow{4}{*}{Q(4,5)}                                       
   & $1.852$&$\mathbf{0.6152}$&$25.2049$         &$0$&$24.6621$         &$1$ \\
   & $1.718$&$4.5963	$        &$\mathbf{25.1177}$&$0$&$24.3975$         &$1$ \\
   & $1.666$&$8.1964	$        &$25.2180$         &$0$&$\mathbf{24.0779}$&$1$ \\
   & $2.0$  &$22.1417$        &$25.4424$         &$0$&$24.9559$         &$1$ \\
   \hline
   
   \multirow{4}{*}{Q(4,6)}                                       
   & $1.831$&$\mathbf{0.3366}$&$25.2581$         &$0$&$24.8755$         &$1$ \\
   & $1.414$&$31.6263$        &$\mathbf{24.2758}$&$0$&$24.2526$         &$1$ \\
   & $1.422$&$30.7519$        &$24.3959$         &$0$&$\mathbf{23.9850}$&$1$ \\
   & $2.0$  &$17.6047$        &$25.3424$         &$0$&$25.1857$         &$1$ \\
   \hline
   
   \multirow{4}{*}{Q(4,7)}          
   & $1.836$&$\mathbf{0.1665}$&$25.2507$         &$0$&$25.0542$         &$1$ \\
   & $1.467$&$25.9829$        &$\mathbf{24.3400}$&$0$&$24.2350$         &$1$ \\
   & $1.460$&$26.7079$        &$24.3655$         &$0$&$\mathbf{24.2174}$&$1$ \\
   & $2.0$  &$20.8774$        &$25.2707$         &$0$&$25.2689$         &$1$ \\
   \hline
   
   \multirow{4}{*}{Q(4,8)}                                       
   & $1.836$&$\mathbf{0.0844}$&$25.2269$         &$0$&$25.0775$         &$1$ \\
   & $1.522$&$20.4675$        &$\mathbf{24.4026}$&$0$&$\mathbf{24.1908}$&$1$ \\
   & $2.0$  &$20.3949$        &$25.2519$         &$0$&$25.2763$         &$1$ \\
   \hline
   
   \multirow{4}{*}{Q(4,9)}                                       
   & $1.834$&$\mathbf{0.0427}$&$25.2008$         &$0$&$25.2364$         &$1$ \\
   & $1.579$&$15.1358$        &$\mathbf{24.4734}$&$0$&$24.4326$         &$1$ \\
   & $1.586$&$14.5189$        &$24.4967$         &$0$&$\mathbf{24.3751}$&$1$ \\
   & $2.0$  &$19.9790$        &$25.2250$         &$0$&$25.2300$         &$1$ \\
   \hline
   
   \multirow{4}{*}{Q(4,10)}                                      
   & $1.834$&$\mathbf{0.0214}$&$25.1843$         &$0$&$25.1848$         &$1$ \\
   & $1.641$&$9.8505	$        &$\mathbf{24.4981}$&$0$&$24.4702$         &$1$ \\
   & $1.643$&$9.6932	$        &$24.5073$         &$0$&$\mathbf{24.4485}$&$1$ \\
   & $2.0$  &$19.6907$        &$25.1871$         &$0$&$25.1280$         &$1$ \\
   \hline
   
   \multirow{4}{*}{Q(5,10)}         
   & $1.915$&$\mathbf{0.0228}$&$25.1996$         &$0$&$25.1956$         &$1$ \\
   & $1.445$&$37.2458$        &$\mathbf{25.1653}$&$0$&$25.1708$         &$1$ \\
   & $1.503$&$31.3220$        &$25.1681$         &$0$&$\mathbf{25.0628}$&$1$ \\
   & $2.0$  &$19.6955$        &$25.1871$         &$0$&$25.1280$         &$1$ \\
   \hline
   
  \end{tabular}
  \begin{tablenotes}
   \item[$\dagger$] Scaling factor, $S = n \implies S = \frac{\max (\text{fp})}{\max(\text{lns})}$
   \item[$\bigstar$] $N_z = $ Number of values rounded to zero in the real domain
  \end{tablenotes}
 \end{threeparttable}
\end{table}

%\begin{figure}
% \centering
% \includegraphics[scale=0.6]{all_mindata.eps}
% \caption{Bases with minimum representation and arithmetic errors for all number formats of 
% Table~\ref{tab:number_format}.}
% \label{fig:all_mindata}
%\end{figure}
%% The plot in Fig.~\ref{fig:lnsquant_scaled_mindata} shows that no one base gives minimum 
%%representation/quantization 
%% error for different word lengths and thus the selection of word length is important in 
%%determining the optimal base.
%
%The best bases corresponding to each type of error is also shown graphically in 
%Fig.~\ref{fig:all_mindata}, which 
%further emphasizes improvements given by smaller bases for low precision, shorter word lengths.

\section{LNS Addition and Subtraction}
\label{sec:lns_addsub}

The set of LNS numbers is not closed under addition and subtraction. For example, if the exponent 
consists only of integer bits (i.e. $f = 0$), then base-$2$ has real-domain numbers $\{1,2,4,8,16 
\dots\}$. The sum of two numbers in this set may not be a member of the set. For example, $1+2 = 
3$, but $3$ is not an element of the set. In practical LNSs, we round the results of addition and 
subtraction to the nearest representable LNS number. Thus, the $+$ and $-$ operations are not exact 
in LNS, leading to rounding error in these fundamental operations. These rounding errors are 
inherent in logarithmic number systems, regardless of how addition and subtraction are implemented.

However, different bases lead to the LNS containing different sets of values. For example, the set 
of LNS numbers with an integer exponent and base-$\sqrt{2}$ contains numbers of the form

\begin{equation*}
\{1, \sqrt(2), 2, 2\sqrt(2), 4 \dots \}.
\end{equation*}

\begin{figure}[h]
	\centering
	\includegraphics[scale=0.5]{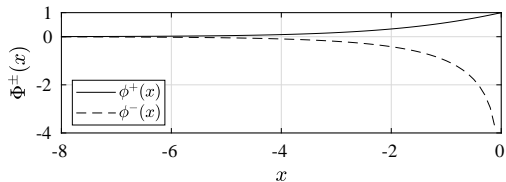}
	\caption{The $\Phi(x)$ function for performing LNS addition and subtraction.}
	\label{fig:lns_addsub_funcplot}
\end{figure}

In this number system, if we add $1+2=3$, the nearest representable LNS number is $2\sqrt(2) 
\approx 2.828$. Thus, both the set of representable values and the size of rounding errors depends 
on the base. A common way to compute LNS addition and subtraction uses the identities 
\cite{kenny_2008}:

\begin{equation}
	m_{\text{add/sub}} = \max{(m_x,m_y)} + \Phi(x)
	\label{eq:log_addsub}
\end{equation}
where
\begin{equation}
	\Phi(x) = 
	\begin{cases}
	\Phi^+(x) = \log_b|1+b^{-x}| \quad S_\Phi = 0  \\
	\Phi^-(x) = \log_b|1-b^{-x}| \quad S_\Phi = 1 \\
	\end{cases}.
	\label{eq:phi}
\end{equation}
where $b$ is the base, $S_\Phi = s_x \oplus s_y \oplus \textit{op}$ (\textit{op} is $0$ for addition and $1$ for 
subtraction) and $x = |m_x-m_y|$ (where $x \le 0$).  The functions $\Phi^+(x)$ and $\Phi^-(x)$ are non-linear, as 
shown in Fig.~\ref{fig:lns_addsub_funcplot}

The challenge of computing function $\Phi(x)$ can offset the gains achieved by other useful 
properties of LNS. A number of methods exist such look-up tables, CORDIC, Taylor series expansion 
or function approximation \cite{Mansour2015J, kenny_2008}. We use the 
look-up table method, which requires $O(2^{i+f})$ space complexity, and is
suitable for word lengths up to $20$-bits \cite{Coleman2008J}. For the
low precision word sizes we consider, in the range $5-16$ bits, the
look-up tables are relatively small.

We represent the $\Phi^+$ and $\Phi^-$ functions as look-up tables for different bases.
% As shown in Equation (\ref{eq:phi}), these tables are indexed by $x = |m_x-m_y|$.
LNS with fixed-point exponent is not closed under addition or subtraction, so we round the
values in the tables to the nearest representable LNS number. Different LNS bases
lead to number systems containing different sets of numbers. Thus, these look-up
tables are different for each base, and therefore the rounding errors are
different for each base.

Typically, in an LNS arithmetic unit, look-up tables are implemented as simple read-only memories 
(ROMs). However, in some cases, the table entries for the two functions are trivial with a very few 
non-zero values and thus realizing them as logic gates can be beneficial than implementing them as 
ROMs.

\subsection{Realization of Large LNS/Small Linear Values in the Subtraction Tables}
\label{subsec:lns_addsub_subtable_cornercases}

According to Equation~\ref{eq:phi}, the first entry in the tables correspond to when $b^{m_x}$ and 
$b^{m_y}$ are equal. This entails two things. First, output of the subtraction of such numbers will 
be zero and that for \phiminus, the first value corresponds to $\log_b|0|$, which is $-\infty$. 
Similar to this, for the next few values, the difference between $b^{m_x}$ and $b^{m_y}$ is very 
small implying that subtraction between these two values results in a very small number. This also 
means that $\log_b|1-b^{-x}|$ of Equation~\ref{eq:log_addsub} returns a very large negative 
value. 

These large negative values in the subtraction table correspond to very small values in the real 
domain, which are smaller than the smallest representable number. This number corresponds to 
$b^{-(2^i-2^{-f})}$. One needs to make a decision about their representation in the tables. 
In the real domain, values that are smaller than the smallest representable number are either 
rounded to zero or the smallest representable number, depending on the distance between the number 
itself and the two adjacent numbers. However, one cannot represent real zero in the LNS domain as 
that is equivalent to $-\infty$. Thus, the underlying question is the representation of linear 
zero/LNS $-\infty$ in the subtraction tables. 

Since in these cases the output of subtraction either results in a zero or a very small value that 
is to be rounded to zero in the real domain, we can utilize this by storing zeros in the table and 
having a small logic circuit that can detect such cases and forcing the output to zero. Storing 
zeros simplifies the logic circuitry synthesized for realizing the subtraction tables.

\subsection{Evaluation of Arithmetic Error}
\label{subsec:lns_addsub_aritherr_eval}

We computed the arithmetic error as a percentage of $\textit{ULP}_{\textit{LNS}}$ in 
the log domain (where $\text{\LULP} = 2^{-f}$) and considered bases ranged from $\sqrt{2}$ and $2$, 
similar to what was done for representation error. The rounding was computed in the log domain with 
the error measured as defined in Section~\ref{subsec:lns_rounding}. 

\begin{figure}
	\centering
	\includegraphics[scale=1.0]{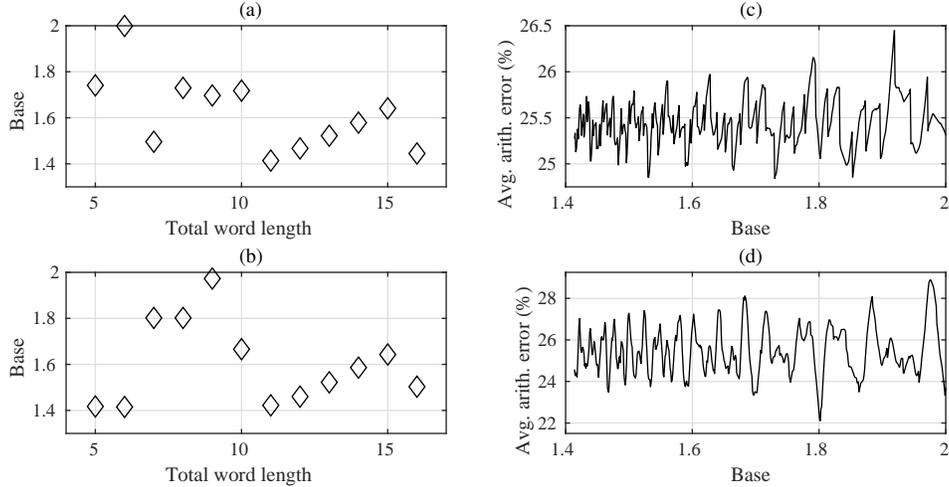}
	\caption{Bases with minimum arithmetic error for various word lengths for (a) addition and (b) subtraction tables. 
	Average absolute arithmetic error for $Q(4,3)$ as a percentage of one \LULP~(in the log domain) for (c) addition and 
	(d) subtraction tables. 
	}
	\label{fig:lnsadd_sub_table_mindata}
\end{figure}

The result is shown in Fig.~\ref{fig:lnsadd_sub_table_mindata}(a,b) for various word lengths 
ranging from $5$ -- $16$ bits (including one sign bit) where only the base that achieves the 
minimum arithmetic error for a given word length is shown. The error measures do not include errors 
generated when values are rounded to zero in the real domain and remain within the bounds given by 
Equation~\ref{eq:mult_error_range}. Also shown in Fig.~\ref{fig:lnsadd_sub_table_mindata}(c,d) is 
the variation in the average absolute arithmetic error as a percentage of \LULP~(in the log domain) 
for $Q(4,3)$ and that for this particular case, base-$1.730$ and base-$1.802$ provide the lowest 
arithmetic error. This is also indicated in Table~\ref{tab:perc_phi_bases}, where the results for 
arithmetic error were presented along with representation error. This is a $2\%$ and $7\%$ 
improvement for the addition and subtraction tables over standard base-$2.0$.

When rounding to nearest, the average arithmetic error will be around $25\%$ of the \LULP~with the 
maximum absolute error being $\le 0.5\times$ \LULP. Figure~\ref{fig:lnsadd_sub_table_mindata}(a,b) 
shows that base-$2$ is not the best choice, especially for smaller word lengths. 

As mentioned in Section~\ref{sec:quant_float}, the errors generated due to representation an input 
distribution is highly domain dependent. However, the values in addition/subtraction tables are 
dependent on the particular arithmetic operation and the base in which they are represented and 
this result shows that rounding error for LNS arithmetic can be reduced by selection of an 
appropriate base.

However, different bases give minimum rounding error for the addition and subtraction tables for 
different word lengths. These bases do not correspond to the bases which give minimum 
representation error as shown in Fig.~\ref{fig:lnsquant_scaled_mindata}. In order to compare the 
arithmetic and representation error for different interesting bases, Table~\ref{tab:perc_phi_bases} 
shows them for each of these bases. It also identifies the best base for each type of error. Also 
shown are the number of values that are rounded to zero in the real domain. These values are not 
counted when computing the arithmetic error, for reasons explained in 
Section~\ref{subsec:lns_rounding}.

%including the first value of the subtraction table which is $-\infty$ and is also represented by a zero.

Typically, for the fixed point exponent, rounding to nearest results in a maximum error of one half 
of the \LULP~and average error of one quarter of the \LULP. However, it can be seen in 
Table~\ref{tab:perc_phi_bases} that for very short word lengths, the average arithmetic error for 
bases other than base-$2$ is much lower than $0.25\times 
\text{\LULP}$, specially for the subtraction table. For e.g., base-$1.417$ gives an arithmetic 
error of $17.1544\%$ of \LULP~for $Q(2,2)$ while base $1.415$ gives $19.9971\%$ of \LULP~(in the 
log domain) for $Q(2,3)$. These errors converge to around $0.25\times\text{\LULP}$ as we increase 
the word length.

Furthermore, as compared to base-$2$, significant improvements can be achieved as compared to 
base-$2$ of up to $17.6\%$ and $48.2\%$ for addition and subtraction tables, respectively for 
$Q(2,2)$. It can also be seen that the representation error for bases which give minimum 
arithmetic error is typically within the $25\%$ range expected of average rounding error, bar only 
one exception when the word length specification matches that of half-precision floating point 
($Q(5,10)$).

\section{Synthesis of Addition and Subtraction Tables}
\label{sec:lns_addsub_synthesis}

A number of methods exist in literature for realizing Equation~\ref{eq:log_addsub} in hardware, as 
mentioned in Section~\ref{sec:lns_addsub}. However, since we are dealing with short word lengths, 
table based approach is the preferred choice \cite{Coleman2008J}.

The functions of \Phipm are evaluated by storing pre-computed values in a look-up table and then 
reading the appropriate value based on the input. The values in the table are quantized to the 
given word length, various choices for which are listed in Table~\ref{tab:number_format}. 
Typically, look-up tables are realized using read-only ROMs. However, if the table entries are 
fairly simple, they can be realized using logic gates rather than ROMs. One can also realize some 
columns of the table using logic gates and others using a ROM depending on the complexity of the 
columns.

Thus, to realize tables, one has three options:

\begin{itemize}
 \item Synthesize the whole table as a ROM
 \item Synthesize the whole table using logic gates
 \item Synthesize parts of the table using logic gates and parts using a ROM
\end{itemize}

A simple, fixed ROM usually consume only one transistor per bit for data storage 
\cite{MemoryMicroProcessorASIC_Book}. However, even a simple fixed ROM will have additional 
circuitry as a row/column decoder, sense amplifiers and precharge transistors \cite{Yang2003J}. A 
simple differential sense amplifier will consist of five transistors \cite{Sedra2019B} while one 
needs one precharge transistor per one column of bits. The address decoding for a memory is 
typically performed by a combination of row and column decoders if memories are arranged as a 
matrix of words rather than a one dimensional array of individual words or only one address decoder 
if arranged as a one dimensional array structure. 

In the case of realizing the tables using logic gates, one can either synthesize the logic for 
individual column of bits independently or combine the columns to enable the synthesis tool to 
share logic to realize various columns of bits. In our experiments, the total synthesis cost, which 
we present in terms of transistors, for a combined synthesis of columns of bits is lower than 
synthesizing individual columns separately. We use the open source \textbf{ABC} synthesis tool 
\cite{abc_reference} and provide the tool with tables described using the Berkeley logic 
interchange format (BLIF) file \cite{Blif_Paper1992} and a library containing only NAND, NOR and 
Inverter gates to enable us to estimate the number of transistors. Each NAND and NOR gate is 
typically made up of four transistors while the inverter is made up of two transistors 
\cite{Sedra2019B}. This can be described as a first degree estimation as we do not take 
into account the differences in the sizes of the two different types of transistors, the NMOS 
(N-type Metal Oxide Semiconductor transistor) and PMOS (P-type metal oxide semiconductor 
transistor) which are used to realize the gates using the complementary MOS technology 
\cite{Sedra2019B}. 

However, the question about when to use logic gates and when to use ROM will depend on a threshold. 
The threshold will be the point when the cost of implementing using gates goes beyond that of ROM 
which has a constant cost of implementation. 

Consider the example of $8$-bit LNS which requires only $128\times 7$ with the sign bit not stored. 
If realized as an array structure, the cost for such a ROM will be $896$ transistors for the 
storage, $7$ transistors for the precharge, $5\times7 = 35$ transistors for sense amplifiers and 
$1022$ transistors for a $7\times 128$ decoder. The implementation cost for the decoder was 
estimated by synthesizing the decoder using the ABC synthesis tool which was provided a BLIF 
file describing the truth table of it. This makes for a total cost of $1960$ transistors. 

However, if we implement the memory as a matrix of $16\times8$ words, one will need a $4\times16$ 
row and $3\times8$ column decoder which cost $122$ and $54$ transistors respectively, making the 
total cost of the memory to be: $896+8\times7 (\text{precharge})+8\times7\times5 (\text{sense 
amplifier}) + 122 + 54 = 1408$ transistors. This number can be further reduced by investigating the 
values each of the addition and subtraction table takes. For the subtraction table, one needs all 
$7$ bits of output, however, for the addition table, one can discard the two most significant bits 
(as they are zero) and only store the least significant five bits. This is also evident by 
Fig.~\ref{fig:lns_addsub_funcplot} where the range of \phiplus is limited between zero and one 
while that of \phiminus is much larger. This can reduce the total cost for the addition table to 
$1056$ transistors. All of this is summarized in Table~\ref{tab:rom_trans}.

\begin{table}
	\centering
	\caption{Implementation cost of various options for a $7\times128$ ROM based table implementation in terms of 
	number of 
	transistors.}
	\label{tab:rom_trans}
	\begin{tabular}{cccccc}
		\hline
		Type of ROM                         & Storage & Precharge & Sense amplifier & Decoder        & Total  \\
		\hline
		ROM with one decoder 				& $ 896 $ & $ 7 $     & $ 35 $          & $ 1022 $       & $1960$ \\
		ROM with two decoders               & $ 896 $ & $ 56 $    & $ 280 $         & $ 122+54=176$  & $1408$ \\
		$5$-bit ROM with two decoders       & $ 640 $ & $ 40 $    & $ 200 $         & $ 122+54=176$  & $1056$ \\
		\hline
	\end{tabular}
\end{table}

\begin{figure}
 \centering
 \includegraphics[scale=0.8]{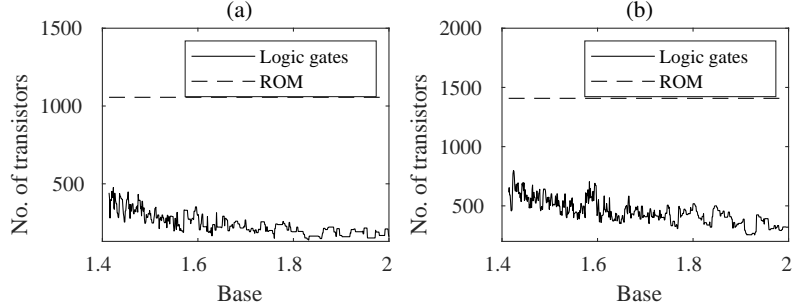}
 \caption{Number of transistors required to realize $8$-bit (a) addition and (b) subtraction tables 
 for all bases 
  between 
  $\sqrt{2}$ and $2$.}
 \label{fig:singlerom_logicgates}
\end{figure}

Comparing the cost of ROM with logic gates for realizing complete table shows that for all bases, 
tables realized using gates has lower area as compared to ROM, as shown in 
Fig.~\ref{fig:singlerom_logicgates} for $8$-bit LNS for all combination of bases between $\sqrt{2}$ 
and $2$, with a minimum savings of $55\%$ and  $43\%$ for the addition and subtraction tables, 
respectively. 

\subsection{Tables using Gates and ROM}
\label{subsec:tables_gates_rom}

Another approach to realize these tables is to realize some columns using logic 
gates and some using ROM. This may appear beneficial because for those columns with only a few 
$1s$, a gate based realization will result in very small circuits while for the columns where 
changes in the value are too frequent, a ROM based implementation may be better or of similar 
complexity as compared to a gate based implementation.

\begin{figure}
	\centering
	\includegraphics[scale=0.8]{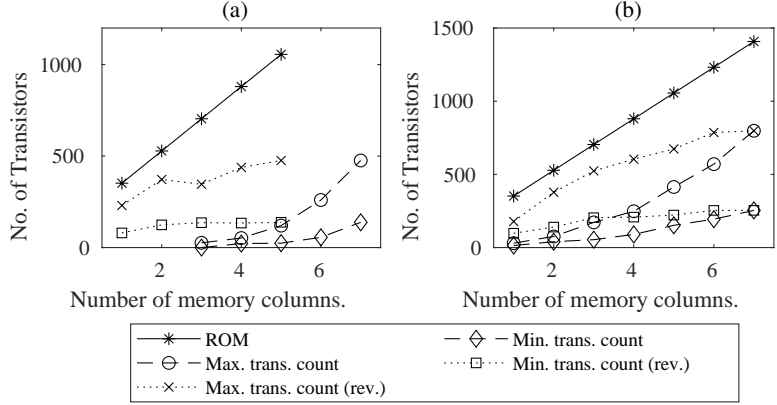}
	\caption{Number of transistors required for realizing successively increasing number of columns for $8$-bit (a) 
	addition 
	and (b) subtraction tables for bases with minimum ($\text{base}_{\phi^+}: 1.832$, $\text{base}_{\phi^-}: 1.913$) 
	and maximum transistor count ($\text{base}_{\phi^+}: 1.423$, $\text{base}_{\phi^-}: 1.424$).}
	\label{fig:singlerom_logicgates_columns}
\end{figure}

However, experimental evaluation shows that is not the case, as shown in 
Fig.~\ref{fig:singlerom_logicgates_columns} for select bases which give the minimum and maximum 
transistor count for each table. The x-axis in each sub-figure shows the number of memory 
columns where \textbf{rev.} means that columns are combined starting from least significant column 
while the other alternative shows the starting point is the most significant columns. For the 
addition table, the two most significant columns are not shown as they are all zeros.

Figure~\ref{fig:singlerom_logicgates_columns} shows that the maximum transistor count for the 
addition table increases from $26$ transistors for a single column to $476$ transistors for all 
columns combined using logic gates for base-$1.423$. This is an improvement of $55\%$ to $93\%$. 
Similarly, the savings in transistor count for the subtraction case, in the worst case, is $43\%$ 
to $91\%$. Even when combining the columns starting from the most significant column results in 
lower transistor count as compared to realizing the tables using a ROM.
 
Thus, considering only the perspective of transistor count, realizing the tables using only logic 
gates results in a lower transistor count as compared to a simple ROM. 
Fig.~\ref{fig:singlerom_logicgates_columns} also indicates that for $8$-bit LNS, a base of $1.832$ 
and $1.913$ gives minimum transistor count for addition and subtraction table, respectively. Thus 
with a selection of base other than base-$2$, one can also reduce the implementation cost.

\subsection{Delay of Logic Circuit for Tables}
\label{subsec:lns_synthesis_dly}

Typically, for a given word length, the latency of reading a word from a ROM will be the same 
across all reads, specially if the word being read are in different rows. In a typical 
implementation, the row address is asserted first to read the whole row (consisting of multiple 
words) into the sense amplifiers from where the exact word is selected using the column address. 
The latency of the two decoders ($4\times 16$ and $3\times 8$, as given in 
Section~\ref{sec:lns_addsub_synthesis}) each for $7$- and $5$-bit memories is $4.7$ and $3.3$ time 
units respectively. Only considering the sum of the delay of the two decoders and ignoring the 
delay associated with the actual reading of the value from transistors, sense amplifiers and any 
refresh times required, the total delay is $8$ time units.

\begin{figure}
	\centering
	\includegraphics[scale=0.8]{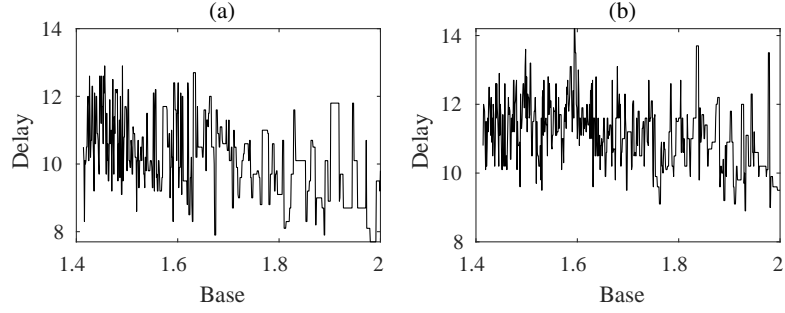}
	\caption{Delay associated with logic circuit required to realize $8$-bit (a) addition and (b) subtraction tables 
	for all 
	bases between $\sqrt{2}$ and $2$.}
	\label{fig:singlerom_logicdelay}
\end{figure}

However, the delays associated with realizing the tables as logic gates for different bases and 
$8$-bit output is shown in Fig.~\ref{fig:singlerom_logicdelay}. The difference between the delay 
associated with logic gates and ROM based implementation is not significant. With a very 
significant difference between the transistors required to realize the two alternatives for 
realizing the addition and subtraction tables, as shown in Figs.~\ref{fig:singlerom_logicgates} 
and \ref{fig:singlerom_logicgates_columns}, it can be concluded that it is better to use logic 
gates for realizing the tables for $8$-bit LNS.

\section{Deciding the Trade-off for Addition and Subtraction Tables}
\label{sec:tradeoff_lnsaddsub}

Four parameters have been discussed in the previous sections. They are:

\begin{itemize}
	\item Representation error
	\item Arithmetic error for \phiplus and \phiminus
	\item Area in terms of transistors count for realizing tables using gates for \phiplus and \phiminus
	\item Delay for realizing tables using gates for \phiplus and \phiminus
\end{itemize}

For each parameter, there are unique bases in the LNS that gives the corresponding minimum value. 
The most important of the these parameters are the two types of errors. 

%Figure~\ref{fig:radar_plot_normprec} shows the trade-offs involved when considering different word 
%lengths. Each point in the figure represents minimum value for the corresponding 
%parameter.
%
%\begin{figure}
%	\centering
%	\includegraphics[scale=0.8]{radar_plot_normprec.eps}
%	\caption{Radar plot showing the trade-offs between various parameters for all word lengths.}
%	\label{fig:radar_plot_normprec}
%\end{figure}

\begin{figure}
	\centering
	\includegraphics[scale=0.8]{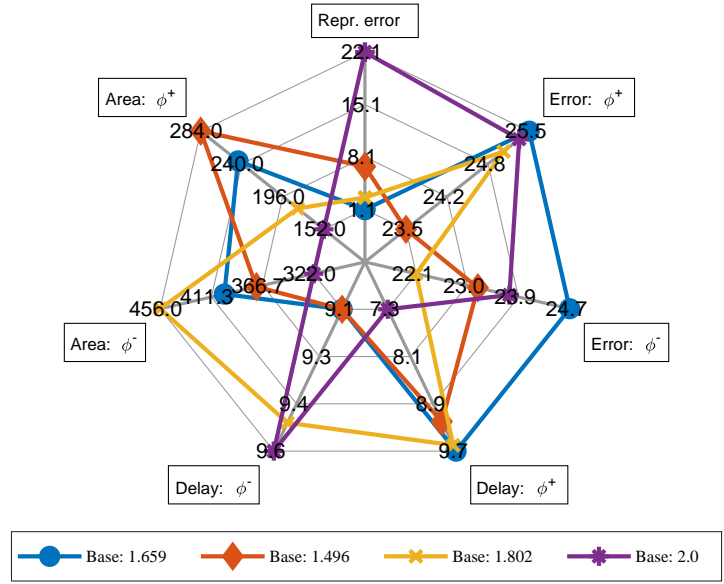}
	\caption{Radar plot showing the trade-offs between various parameters for Q(3,3).}
	\label{fig:radar_plot_normprec_q3_3}
\end{figure}

The trade-offs involved for one particular word length of $Q(3,3)$ is shown in  
Fig.~\ref{fig:radar_plot_normprec_q3_3} for various bases. Similar trade-offs exist for other word 
lengths. All parameters have their importance in the overall design and implementation and each 
base is optimal in one of the four parameters defined above but none of the base is optimal in all 
of the four parameters. This figure emphasize our point of the need of analyzing different bases 
and the priorities while designing and implementing any system. However, it has 
been clearly shown that base-$2.0$ is not optimal under different word lengths and for different 
parameters.

\section{Higher Precision Tables for Very Short Word Lengths}
\label{sec:high_precision_outputs}

As mentioned earlier, a number of works use very short word lengths for realizing different systems 
for deployment on edge devices. However, short word lengths can introduce significant precision 
errors for addition and subtraction. One way to circumvent this is to use tables with higher 
precision for table values (referred from here on as mixed precision tables for exploratory 
convenience) and rounding the result only after all required operations for, say, computing the sum 
of products of data and coefficients in an FIR filter is complete. Thus, there is only one source 
of arithmetic error at the output (at the cost of higher hardware cost).%computing the feature maps 
%of a single layer of a 
%neural network or the output of a digital filter.

\begin{figure}
	\centering
	\includegraphics[scale=0.8]{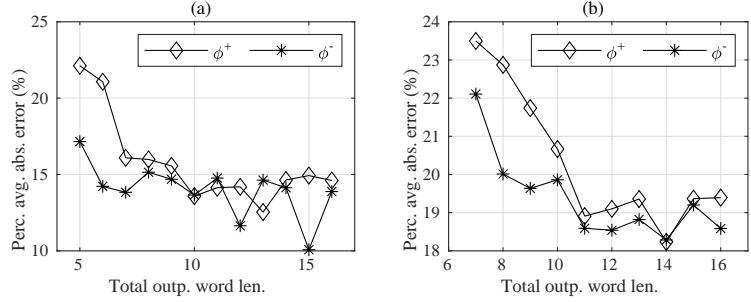}
	\caption{Best case average absolute arithmetic error as a percentage of the ULP~when the addition and subtraction 
	mixed precision tables are quantized to increasing fractional word length for (a) Q(2,2) and (b) Q(3,3) number formats 
	identified in Table~\ref{tab:number_format}.}
	\label{fig:lnsadd_sub_mprec_table_mindata}
\end{figure}

The number format given in Table~\ref{tab:number_format} was used for the analysis with fractional 
word length successively increased up to a maximum of total word length of $16$ (including the sign 
bit) for each case. For e.g., for the case of $Q(2,2)$, the input word length is 
$5$-bits while the output word length is varied from $5$--$16$ bits. The results for two different 
total input word length, Q(2,2) and Q(3,3), are shown in 
Fig.~\ref{fig:lnsadd_sub_mprec_table_mindata}, where each value represents the minimum percentage 
average absolute arithmetic error. For each value, the base that gives this minimum value is 
unique. As discussed earlier in Section~\ref{subsec:lns_addsub_aritherr_eval}, for shorter word 
lengths, the arithmetic error can be reduced much lower than $25\%$ of \LULP~for bases other than 
base-$2$.

\begin{figure}
	\centering
	\includegraphics[scale=0.8]{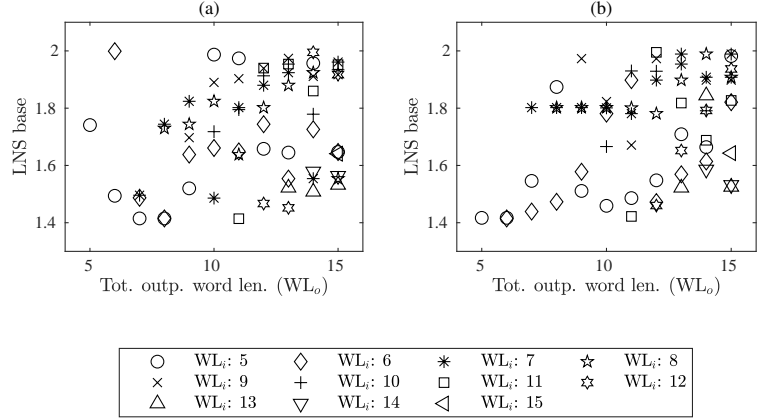}
	\caption{Bases with minimum arithmetic error for various word lengths for (a) addition and (b) 
	subtraction mixed precision tables.}
	\label{fig:bestbase_mprec_aritherror}
\end{figure}

For each number format given in Table~\ref{tab:number_format} and all output word 
lengths, there is no single optimal base which gives the least arithmetic error for various 
combinations of integer and fractional bits. This phenomenon is further emphasized in 
Fig.~\ref{fig:bestbase_mprec_aritherror} which identifies the various bases which give the least 
arithmetic error for total output word length. The legend shows the original total word length 
(without extended fractional bits). For e.g., the plot shows the optimal base for output word 
lengths ranging from $5$--$16$ bits for the case of $Q(2,2)$. Again it re-emphasizes the earlier 
reached conclusion that bases other than base-$2$ should be considered when designing a system as 
it can give lower error and implementation cost.

\subsection{Synthesis of Higher Precision Tables for Very Short Word Lengths}
\label{subsec:high_precision_outputs_synthesis}

The other key issue to deal with extra fractional bits is the size of these tables, once they are 
realized in hardware. As shown in Figs.~\ref{fig:singlerom_logicgates} and 
\ref{fig:singlerom_logicgates_columns}, it is better to realize the tables using gates. However, 
with increasing output word length for the tables, one needs to analyze whether the conclusion 
reached in Section~\ref{sec:lns_addsub_synthesis} still stands. For this, the addition and 
subtraction tables were synthesized using both ROM and standard cell gates for all input and output 
word length combinations. Table~\ref{tab:rom_mprec_trans} shows the number of transistors required 
to realize the ROMs based on the model earlier defined in Section~\ref{sec:lns_addsub_synthesis}. 
The size of decoders for the various tables were chosen so as to achieve minimum transistor count 
for their implementation.

\begin{table}
	\centering
	\footnotesize
	\caption{Implementation cost for ROM based table implementation, for mixed precision tables, in terms of number of 
	transistors.}
	\label{tab:rom_mprec_trans}
	\begin{tabular}{C{1.5cm}|C{0.4cm}|C{0.4cm}|C{0.4cm}|C{0.575cm}|C{0.575cm}|C{0.575cm}|C{0.7cm}|
   C{0.7cm}|C{0.7cm}|C{0.8cm}|C{0.8cm}|C{0.80cm}}
		\hline
		\multirow{2}{1.5cm}{Input word length ($\text{WL}_i$)} & \multicolumn{12}{c}{Transistor count for ROM 
		(with various output word lengths ($\text{WL}_o$))} \\
		\cline{2-13}
		     & $5$ & $6$ & $7$ & $8$  & $9$  & $10$ & $11$  & $12$  & $13$  & $14$   & $15$   & $16$   \\\hline
		$5 $ & 168 & 196 & 224 & 252  & 280  & 308  & 336   & 364   & 392   & 420    & 448    & 476    \\
		$6 $ &  -- & 345 & 401 & 457  & 513  & 569  & 625   & 681   & 737   & 793    & 849    & 905    \\
		$7 $ &  -- & --  & 661 & 749  & 837  & 925  & 1013  & 1101  & 1189  & 1277   & 1365   & 1453   \\ 
		$8 $ &  -- & --  & --  & 1408 & 1584 & 1760 & 1936  & 2112  & 2288  & 2464   & 2640   & 2816   \\ 
		$9 $ &  -- & --  & --  & --   & 2736 & 3040 & 3344  & 3648  & 3952  & 4256   & 4560   & 4864   \\ 
		$10$ &  -- & --  & --  & --   & --   & 5844 & 6452  & 7060  & 7668  & 8276   & 8884   & 9492   \\
		$11$ &  -- & --  & --  & --   & --   & --   & 11832 & 12952 & 14072 & 15192  & 16312  & 17432  \\ 
		$12$ &  -- & --  & --  & --   & --   & --   & --    & 25400 & 27640 & 29880  & 32120  & 34360  \\ 
		$13$ &  -- & --  & --  & --   & --   & --   & --    & --    & 52728 & 57016  & 61304  & 65592  \\ 
		$14$ &  -- & --  & --  & --   & --   & --   & --    & --    & --    & 113020 & 121596 & 130172 \\ 
		$15$ &  -- & --  & --  & --   & --   & --   & --    & --    & --    & --     & 237312 & 254080 \\ 
		$16$ &  -- & --  & --  & --   & --   & --   & --    & --    & --    & --     & --     & 506112 \\
		\hline
	\end{tabular}
\end{table}

\begin{figure}
	\centering
	\includegraphics[scale=0.8]{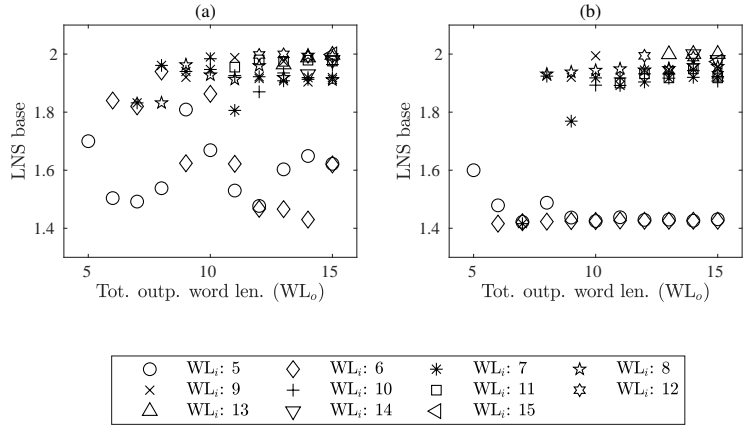}
	\caption{Bases with minimum transistor count for various word lengths for the (a) addition and (b) 
	subtraction mixed precision tables.}
	\label{fig:bestbase_mprec_transcount}
\end{figure}

The addition (\phiplus) and subtraction (\phiminus) tables were also synthesized for all input and 
output word length combinations using standard cell logic gates. For each case, the base which 
gives the minimum transistor count is unique and a distribution of these bases is shown in 
Fig.~\ref{fig:bestbase_mprec_transcount}. This figure, together with the results shown in 
Fig.~\ref{fig:bestbase_mprec_aritherror}, identifies the various bases which give minimum 
arithmetic error and implementation cost, respectively, with lower the word length, the more 
smaller bases give minimum transistor count.

However, arithmetic error is more important than transistor count to ensure correctness of result. 
Although the hardware cost of realizing tables with bases giving minimum arithmetic error is larger 
than the best case in terms of transistor count, it is still significantly better than realizing it 
in ROMs. Results for the subtraction table are shown in 
Table~\ref{tab:phim_gates_mprec_trans_bestbase} as typically it is more complex to realize than the 
addition table and shows an order of magnitude improvement, with a minimum reduction of $168\%$ in 
transistor count to a maximum of nearly $10\times$ reduction.

\begin{table}
	\centering
	\footnotesize
	\caption{Implementation cost for gates based table implementation, for mixed precision table for \phiminus, in 
	terms of number of transistors for bases that give the minimum arithmetic error.}
	\label{tab:phim_gates_mprec_trans_bestbase}
	\begin{tabular}{C{1.5cm}|C{0.4cm}|C{0.4cm}|C{0.4cm}|C{0.575cm}|C{0.575cm}|C{0.575cm}|C{0.7cm}|
   C{0.7cm}|C{0.7cm}|C{0.8cm}|C{0.8cm}|C{0.80cm}}
		\hline
		\multirow{2}{1.5cm}{Input word length ($\text{WL}_i$)} & \multicolumn{12}{c}{Transistor Count for 
		ROM (with various output word lengths ($\text{WL}_o$))} \\
		\cline{2-13}
		& $5$ & $6$ & $7$  & $8$  & $9$  & $10$  & $11$  & $12$  & $13$   & $14$   & $15$  & $16$ \\\hline
		$5 $ & 64 & 96  & 122 & 120 & 156 & 186  & 208  & 214   & 244   & 254   & 282   & 290\\
		$6 $ & -- & 158 & 244 & 264 & 394 & 416  & 486  & 452   & 530   & 614   & 662   & 788\\
		$7 $ & -- & --  & 456 & 454 & 686 & 808  & 834  & 810   & 1114  & 1026  & 1336  & 1484\\
		$8 $ & -- & --  & --  & 494 & 530 & 608  & 872  & 1218  & 1074  & 1768  & 1428  & 2582\\
		$9 $ & -- & --  & --  & --  & 878 & 1326 & 1744 & 2996  & 2122  & 3962  & 3018  & 4460\\
		$10$ & -- & --  & --  & --  & --  & 1746 & 3604 & 4196  & 5150  & 6040  & 4310  & 7528\\
		$11$ & -- & --  & --  & --  & --  & --   & 6132 & 7130  & 8916  & 10080 & 10664 & 12458\\
		$12$ & -- & --  & --  & --  & --  & --   & --   & 10448 & 13072 & 15884 & 16184 & 19132\\
		$13$ & -- & --  & --  & --  & --  & --   & --   & --    & 18962 & 21918 & 27644 & 31194\\
		$14$ & -- & --  & --  & --  & --  & --   & --   & --    & --    & 33438 & 41068 & 50580\\
		$15$ & -- & --  & --  & --  & --  & --   & --   & --    & --    & --    & 63140 & 75726\\
		$16$ & -- & --  & --  & --  & --  & --   & --   & --    & --    & --    & --    & 66608\\
		\hline                                  
	\end{tabular}
\end{table}

It can also be seen from Table~\ref{tab:phim_gates_mprec_trans_bestbase} that as the output word 
length increases, for each input word length, the transistor count does not increase significantly, 
making the mixed precision approach an attractive one given the reduced error measurement. 

\section{LNS Adder and Multiplier Design in Fixed and Floating Point Number System}
\label{sec:lns_v_fixedfloat}

As identified in Sections~\ref{sec:intro} and \ref{sec:related_work}, logarithmic number system 
(LNS) finds application in signal processing and neural networks. Both of these systems are based 
heavily on the operation of addition and multiplication. The key advantage of LNS over fixed and 
floating point stem from its inherent simplification of multiplication into addition. However, the 
drawback with LNS is that it can only approximate the addition operation since addition is not 
closed under LNS, as outlined in Section~\ref{sec:lns_addsub}. 

Since multiplication in LNS is reduced to addition, the circuit to realize an LNS multiplier will 
be similar in complexity to that of a fixed point adder \cite{Kouretas2018C}. This is because the 
exponent of the LNS is typically represented using fixed point. In order to ascertain 
which number system to use from a hardware complexity purpose it is important to realize the 
hardware circuit for LNS addition/subtraction and compare it against fixed point multiplier. 
Since in applications like neural networks the dominant number system used is floating point, 
it is important to realize floating point addition and multiplication circuits to analyze which 
number system provides the smallest circuit. 

\begin{figure}
 \centering
 \includegraphics[scale=0.35]{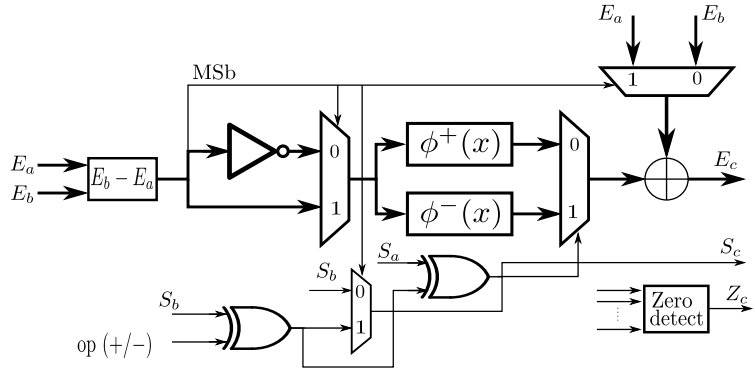}
 \caption{Circuit schematic for LNS addition and subtraction.}
 \label{fig:lns_addsub}
\end{figure}

The circuit to perform addition and subtraction in LNS is given in Fig.~\ref{fig:lns_addsub} and is 
a modified form of the circuit given by Kenny et al in \cite{kenny_2008}. In the circuit, the 
functions of \phiplus and \phiminus are implemented using look-up tables, as mentioned in 
Section~\ref{sec:lns_addsub_synthesis}. The zero detect module shown takes into account whether the 
two numbers being added are both zero (shown by the $z_x$ bit of Equation(\ref{eq:log})), the 
difference of two number is zero or so small that the result will be rounded to zero in the 
real/linear domain, as discussed in Section~\ref{subsec:lns_addsub_subtable_cornercases}.

Addition and multiplier circuits were synthesized for logarithmic, fixed and floating 
point number systems for word length ranging from five to $16$-bits using Cadence Genus with the 
45nm GSCL liberty cell library from the Free45PDK. For LNS, the bases were varied from $\sqrt{2}$ 
to $2$, similar to how LNS is analyzed in the previous sections. The results are shown in 
Table~\ref{tab:lns_v_fixedfloat} where $Q(i,f)$ indicates $i$ integer and $f$ fractional bits and 
does not include the sign bit. For fixed point, it needs to be read as $Q(i+1,f)$ since $2's$ 
complement number representation is assumed and for floating point, results for both adder 
(Float.A) and multiplier (Float.M) are presented.

As evident in Table~\ref{tab:lns_v_fixedfloat} the LNS adder requires lower or similar area as 
compared to fixed point multipliers for short word lengths up to $Q(4,4)$. Although the 
floating point multiplier is the smaller of all circuits (for most of the cases), the LNS 
multiplier and fixed-point adder will be more simpler and smaller than it, thus overall, floating 
point will not have the smallest of the circuit. Also important to realize is that for LNS, the 
same circuit is used to perform both addition and subtraction which opens up the possibility of 
re-using the same circuit for the two operations if the input sample rate allows for some level of 
multiplexing. Furthermore, for LNS, the adder and subtractor were realized using the look-up table 
approach for implementing the functions \phiplus and \phiminus in Equation (\ref{eq:log_addsub}). 
They can also be realized using different techniques like weighted sum of bit-products 
\cite{kenny_2008} which can reduce the area complexity. However, in this work all error 
calculations are based on table based implementations, the same method was used for realizing 
circuits for consistency.

\begin{table}
    \centering
    \footnotesize
%    \begin{threeparttable}
    \caption{A comparison of synthesis results for fixed point multiplier, floating point adder and 
    multiplier and LNS adder circuits.}
    \label{tab:lns_v_fixedfloat}
    \begin{tabular}{C{0.95cm}|l|C{0.75cm}C{0.75cm}C{0.75cm}|C{0.95cm}|l|C{0.9cm}C{0.75cm}C{0.75cm}}
     \hline
     Number format & Circuit type & Area ($\mu m^2$) & Delay ($ps$) & Power ($\mu W$) & Number 
     format 
     & Circuit type & Area ($\mu m^2$) & Delay ($ps$) & Power ($\mu W$)\\
        \hline
        \multirow{3}{*}{$Q(2,2)$ } & Fixed         & 495   & 874  & 101 & 
        \multirow{3}{*}{$Q(2,3)$ } & Fixed         & 688   & 1035 & 146 \\
                                   & Float.M       & 415   & 991  & 80  & 
                                   & Float.M       & 576   & 1467 & 101 \\
                                   & Float.A       & 565   & 1868 & 108 & 
                                   & Float.A       & 706   & 2325 & 120 \\
                                   & LNS ($2.0$)   & 445   & 1316 & 109 & 
                                   & LNS ($2.0$)   & 577   & 1650 & 148 \\
                                   & LNS ($1.421$) & 423   & 1214 & 107 & 
                                   & LNS ($1.443$) & 564   & 1583 & 147 \\\hline
        \multirow{3}{*}{$Q(3,3)$ } & Fixed         & 890   & 1190 & 195 & 
        \multirow{3}{*}{$Q(4,3)$ } & Fixed         & 1143  & 1328 & 251 \\
                                   & Float.M       & 680   & 1504 & 126 & 
                                   & Float.M       & 730   & 1620 & 135 \\
                                   & Float.A       & 913   & 2922 & 170 & 
                                   & Float.A       & 1002  & 3314 & 170 \\
                                   & LNS ($2.0$)   & 744   & 1900 & 170 & 
                                   & LNS ($2.0$)   & 883   & 2170 & 197 \\
                                   & LNS ($1.718$) & 713   & 1832 & 175 & 
                                   & LNS ($1.999$) & 857   & 2061 & 191 \\\hline
        \multirow{3}{*}{$Q(4,4)$ } & Fixed         & 1332  & 1515 & 303 & 
        \multirow{3}{*}{$Q(4,5)$ } & Fixed         & 1599  & 1753 & 379 \\
                                   & Float.M       & 962   & 1931 & 180 & 
                                   & Float.M       & 1102  & 2307 & 216 \\
                                   & Float.A       & 1145  & 3179 & 222 & 
                                   & Float.A       & 1305  & 3557 & 249 \\
                                   & LNS ($2.0$)   & 1325  & 2477 & 263 & 
                                   & LNS ($2.0$)   & 1622  & 2477 & 263 \\
                                   & LNS ($1.96$)  & 1244  & 2660 & 264 & 
                                   & LNS ($1.948$) & 1572  & 2930 & 347 \\\hline
        \multirow{3}{*}{$Q(4,6)$ } & Fixed         & 1856  & 1890 & 458 & 
        \multirow{3}{*}{$Q(4,7)$ } & Fixed         & 2208  & 2109 & 541 \\
                                   & Float.M       & 1328  & 2518 & 258 & 
                                   & Float.M       & 1553  & 2831 & 318 \\
                                   & Float.A       & 1460  & 4118 & 288 & 
                                   & Float.A       & 1653  & 4688 & 350 \\
                                   & LNS ($2.0$)   & 2643  & 3374 & 552 & 
                                   & LNS ($2.0$)   & 4319  & 4342 & 824 \\
                                   & LNS ($1.92$)  & 2523  & 3316 & 550 & 
                                   & LNS ($1.997$) & 4113  & 4022 & 804 \\\hline
        \multirow{3}{*}{$Q(4,8)$ } & Fixed         & 2531  & 2502 & 646 & 
        \multirow{3}{*}{$Q(4,9)$ } & Fixed         & 2842  & 2652 & 726 \\
                                   & Float.M       & 1796  & 3008 & 380 & 
                                   & Float.M       & 2062  & 3364 & 438 \\
                                   & Float.A       & 1769  & 4528 & 355 & 
                                   & Float.A       & 1958  & 5358 & 423 \\
                                   & LNS ($2.0$)   & 7468  & 5419 & 1269& 
                                   & LNS ($2.0$)   & 13155 & 6328 & 1897\\
                                   & LNS ($1.998$) & 7289  & 5176 & 1234& 
                                   & LNS ($1.975$) & 13032 & 6341 & 1906\\\hline
%        \multirow{3}{*}{$Q(4,10)$} & Fixed         & 3206  & 2926 & 835 & 
%        \multirow{3}{*}{$Q(5,10)$} & Fixed         & 3519  & 2956 & 926 \\
%                                   & Float.M       & 2429  & 3622 & 570 & 
%                                   & Float.M       & 2475  & 3627 & 545 \\
%                                   & Float.A       & 2080  & 4953 & 437 & 
%                                   & Float.A       & 2228  & 5152 & 437 \\
%                                   & LNS ($2.0$)   & 23435 & 7315 & 2829 &
%                                   & LNS ($2.0$)   & 24712 & 6629 & 2187\\
%                                   & LNS ($1.999$) & 23376 & 7252 & 2724 &
%                                   & LNS ($1.999$) & 24547 & 6798 & 2190\\\hline
    \end{tabular}
%    \begin{tablenotes}
%     \item[$\dagger$] Unit of measurement: $\mu m^2$
%     \item[$\bigstar$] Unit of measurement: $ps$
%     \item[$\ddagger$] Unit of measurement: $\mu W$
%    \end{tablenotes}
%    \end{threeparttable}
%    \vspace*{-0.675cm}
\end{table}

\begin{figure}[t]
 \centering
 \includegraphics[scale=0.75]{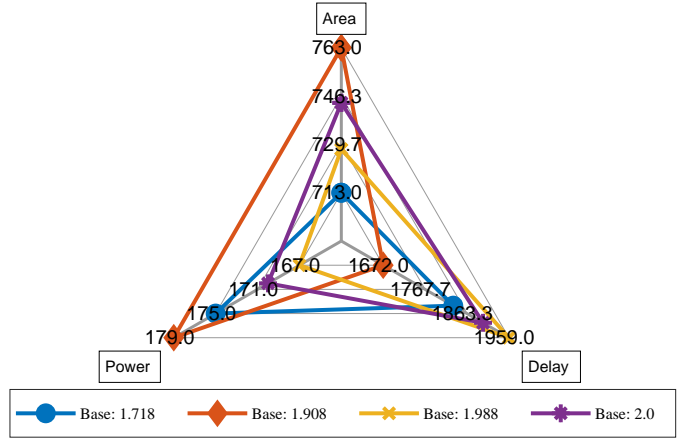}
 \vspace*{-0.75cm}
 \caption{Radar plot showing the trade-offs between the three design parameters for $Q(3,3)$.}
 \label{fig:radar_plot_normprec_q3_3_addsub}
\end{figure}

For each number format, the result for LNS other than base-$2$ is the base which gives the lowest 
area. This does not correspond to the shortest delay through the circuit or the lowest power 
consumption. Fig.~\ref{fig:radar_plot_normprec_q3_3_addsub} gives a better insight into the 
inherent trade-offs involved here.

\section{Finite-Impulse Response (FIR) Filter: An Application}
\label{sec:fir_filter}

\begin{figure}
 \centering
 \includegraphics[scale=1.0]{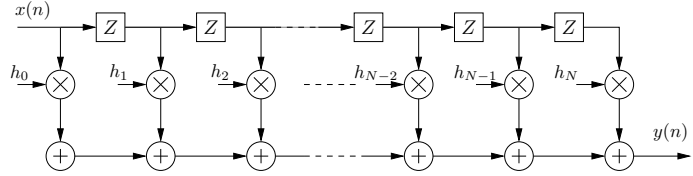}
 \caption{Direct form FIR filter.}
 \label{fig:DF_FIR}
\end{figure}

Frequency selective finite-length impulse response (FIR) filters are used in many applications and 
are therefore a key building block in a digital signal processing system\cite{mitra}. The 
difference equation that defines the output of an $N$th-order (length $N+1$) FIR filter is:
\begin{equation}
y(n) = \sum_{k=0}^{N}h(k)x(n-k)
\label{eq:fir}
\end{equation}
where $y(n)$ is the output sequence, $x(n)$ is the input sequence, and $h(k)$ are the coefficients. 
The direct form of FIR filter to realize Equation(\ref{eq:fir}) is shown in Fig.~\ref{fig:DF_FIR}.

The arithmetic complexity of FIR filters primarily depends on the number of multiplications, which 
according to (\ref{eq:fir}) is proportional to the filter order. The filter order is, for 
linear-phase single stage FIR (SSF) filters, roughly proportional to the inverse of the width of 
the transition band \cite{mitra}. 

%It should be noted that for linear-phase FIR filters, the filter 
%coefficient are symmetric (or anti-symmetric), and, hence, the number of multiplications is 
%$\left\lceil{\frac{N+1}{2}}\right\rceil$.

With LNS reducing multiplication to an addition, it finds application in applications like FIR 
filters where the main operator contributing towards overall hardware complexity is multiplication. 
However, as stated earlier, addition/subtraction in LNS is non-trivial and requires a table based 
implementation as shown in Fig.~\ref{fig:lns_addsub}. As shown in 
Section~\ref{sec:lns_v_fixedfloat}, LNS provides lower hardware cost in terms of area for short 
word length adders/subtractors as compared to its counterparts, it is important to analyze 
whether these benefits extend to FIR filters.

For this analysis, a $11^{\text{th}}$ order FIR filter was implemented in hardware using the same 
setup as described in Section~\ref{sec:lns_v_fixedfloat} with
fixed, floating and logarithmic number systems. The input, output and filter coefficients were 
assumed to be in respective number systems. For LNS, the FIR filter was realized for a select 
number of bases identified earlier in Tables~\ref{tab:perc_phi_bases} and 
\ref{tab:lns_v_fixedfloat}. The results of this analysis is shown in Table~\ref{tab:fir_results}.

The results shown in Table~\ref{tab:fir_results} are consistent with those shown in 
Table~\ref{tab:lns_v_fixedfloat}. FIR filters realized using logarithmic number systems achieve a 
lower area (minimum savings of $2\%$) while also having the shortest delay (minimum improvement of 
$4\%$)and lower power consumption (minimum improvement of $6\%$)as compared to its 
fixed and floating point counterparts up to the word length corresponding to $Q(4,4)$. Overall, LNS 
based FIR filters also are the fastest even though occupies an area that is, at a maximum, nearly 
$3\times$ more than others. 

\begin{table}
 \centering
 \footnotesize
 %\begin{threeparttable}
 \caption{A comparison of synthesis results for FIR filters implemented using fixed, floating and 
 logarithmic number systems.}
 \label{tab:fir_results}
 \begin{tabular}{C{0.95cm}|l|C{0.75cm}C{0.75cm}C{0.75cm}|C{0.95cm}|l|C{0.9cm}C{0.75cm}C{0.75cm}}
  \hline
  Number format & Circuit type & Area ($\mu m^2$) & Delay ($ps$) & Power ($\mu W$) & Number format 
  & Circuit type & Area ($\mu m^2$) & Delay ($ps$) & Power ($\mu W$)\\
  \hline
  \multirow{3}{*}{$Q(2,2)$ } & Fixed         & 4914  & 193 & 1700 & 
  \multirow{3}{*}{$Q(2,3)$ } & Fixed         & 6841  & 197 & 2424 \\
                             & Float         & 7891  & 414 & 3346 & 
                             & Float         & 10877 & 411 & 5200 \\
                             & LNS ($2.0$)   & 3663  & 185 & 1542 & 
                             & LNS ($2.0$)   & 5301  & 189 & 2443 \\
                             & LNS ($1.414$) & \textbf{3454}  & \textbf{185} & \textbf{1297} & 
                             & LNS ($1.417$) & \textbf{4620}  & \textbf{185} & \textbf{1828} 
                             \\\hline
  \multirow{3}{*}{$Q(3,3)$ } & Fixed         & 8964  & 202 & 3239 & 
  \multirow{3}{*}{$Q(4,3)$ } & Fixed         & 11456 & 208 & 4251 \\
                             & Float         & 12554 & 458 & 5694 & 
                             & Float         & 13693 & 524 & 5275 \\
                             & LNS ($2.0$)   & 6691  & 194 & 3007 & 
                             & LNS ($2.0$)   & 8393  & 195 & 3231 \\
                             & LNS ($1.975$) & \textbf{6568}  & \textbf{194} & \textbf{2957} & 
                             & LNS ($1.96$)  & \textbf{8044}  & \textbf{195} & \textbf{3117} 
                             \\\hline
  \multirow{3}{*}{$Q(4,4)$ } & Fixed         & 14104 & 207 & 5269 & 
  \multirow{3}{*}{$Q(4,5)$ } & Fixed         & 16857 & 216 & \textbf{6601} \\
                             & Float         & 16418 & 540 & 6042 & 
                             & Float         & 19250 & 560 & 7907 \\
                             & LNS ($2.0$)   & \textbf{11945} & \textbf{190} & \textbf{4810} & 
                             & LNS ($2.0$)   & 16954 & 195 & 7406 \\
                             & LNS ($1.96$)  & 12436 & 194 & 4969 & 
                             & LNS ($1.998$) & \textbf{16554} & \textbf{195} & 7223 \\\hline
  \multirow{3}{*}{$Q(4,6)$ } & Fixed         & \textbf{20860} & 246 & \textbf{8091} & 
  \multirow{3}{*}{$Q(4,7)$ } & Fixed         & \textbf{23762} & 244 & \textbf{9373} \\
                             & Float         & 23351 & 587 & 9999 & 
                             & Float         & 27019 & 619 & 11221 \\
                             & LNS ($2.0$)   & 27674 & 190 & 12295 & 
                             & LNS ($2.0$)   & 45148 & 195 & 18193 \\
                             & LNS ($1.92$)  & 27183 & \textbf{190} & 12019 & 
                             & LNS ($1.997$) & 44707 & \textbf{190} & 11896 \\\hline
  \multirow{3}{*}{$Q(4,8)$ } & Fixed         & \textbf{27866} & 256 & \textbf{11087} & 
  \multirow{3}{*}{$Q(4,9)$ } & Fixed         & \textbf{31347} & 251 & \textbf{12698} \\
                             & Float         & 30055 & 674 & 13238 & 
                             & Float         & 37063 & 659 & 15873 \\
                             & LNS ($2.0$)   & 81810 & \textbf{189} & 32375 & 
                             & LNS ($2.0$)   & 146520& \textbf{189} & 47353 \\
                             & LNS ($1.998$) & 79898 & 195 & 31760 & 
                             & LNS ($1.975$) & 149660& 185 & 48731 \\\hline
%  \multirow{3}{*}{$Q(4,10)$} & Fixed         & \textbf{35974} & 268 & \textbf{14698} & 
%  \multirow{3}{*}{$Q(5,10)$} & Fixed         & \textbf{39534} & 266 & \textbf{16320} \\
%                             & Float         & 40190 & 659 & 18337& 
%                             & Float         & 41716 & 659 & 17122\\
%                             & LNS ($2.0$)   & 258878& 265 & 85207&
%                             & LNS ($2.0$)   & 261620& \textbf{254} & 76398\\
%                             & LNS ($1.998$) & 257308& \textbf{244} & 88202&
%                             &--             & --    & --  & --\\\hline
 \end{tabular}
% \begin{tablenotes}
%  \item[$\dagger$] Unit of measurement: $\mu m^2$
%  \item[$\bigstar$] Unit of measurement: $ps$
%  \item[$\ddagger$] Unit of measurement: $\mu W$
% \end{tablenotes}
% \end{threeparttable}
%\vspace*{-0.65cm}
\end{table}

The large area for LNS based FIR filters for larger word lengths is due to the 
exponential increase in the memory size required for implementing the \phiplus and \phiminus 
functions, as given in Equation~\ref{eq:log_addsub}. This can be reduced by using approaches like 
the co-transformation approach by Arnold et al. in \cite{Arnold1998J}, the multi-table approach 
proposed by Coleman et al. in \cite{Coleman2008J} or the multipartite table methods proposed by 
Dinechin et al. in \cite{Dinechin2001C}.

\begin{table}
 \centering
 \caption{Cost of conversion between logarithmic and fixed point number systems in terms of area 
 ($\mu m^2$)}
 \label{tab:conv_ckt}
 \begin{tabular}{l|cccccccccc}
  \hline
  & \multicolumn{2}{c}{$Q(2,2)$} & \multicolumn{2}{c}{$Q(2,3)$} & \multicolumn{2}{c}{$Q(3,3)$} & 
  \multicolumn{2}{c}{$Q(4,3)$} & \multicolumn{2}{c}{$Q(4,4)$}\\\hline
  Base &$ 2.0 $&$ 1.414 $&$ 2.0 $&$ 1.417 $&$ 2.0 $&$ 1.975 $&$ 2.0 $&$ 1.96 $&$ 2.0 $&$ 1.96 
  $\\\hline
  Fixed $\rightarrow$ LNS &$ 52  $&$ 45 $&$ 149 $&$ 120 $&$ 310 $&$ 358 $&$ 512 $&$ 609 $&$ 1325 
  $&$ 1429 $ \\
  LNS $\rightarrow$ Fixed &$ 29 $&$ 30 $&$ 68 $&$ 111 $&$ 146 $&$ 185 $&$ 308 $&$ 217 $&$ 608 $&$ 
  651 $ \\\hline
  & \multicolumn{2}{c}{$Q(4,5)$} & \multicolumn{2}{c}{$Q(4,6)$} & \multicolumn{2}{c}{$Q(4,7)$} & 
  \multicolumn{2}{c}{$Q(4,8)$} & \multicolumn{2}{c}{$Q(4,9)$} \\\hline
  Base &$ 2.0 $&$ 1.998 $&$ 2.0 $&$ 1.92 $&$ 2.0 $&$ 1.997 $&$ 2.0 $&$ 1.998 $&$ 2.0 $&$ 1.975 $ 
  \\\hline
  Fixed $\rightarrow$ LNS &$ 2721 $&$ 2715 $&$ 5362 $&$ 5433 $&$ 10879 $&$ 10751 $&$ 19631 $&$ 
  20418 $&$ 38118 $&$ 38783 $ \\
  LNS $\rightarrow$ Fixed &$ 1413 $&$ 1453 $&$ 2609 $&$ 2782 $&$ 5353 $&$ 5265 $&$ 10195 $&$ 10347 
  $&$ 20003 $&$ 20577 $ 
  \\\hline  
 \end{tabular}
\end{table}

Typically conversion can be carried out offline and all calculations can be done online in LNS. 
However, to show the cost of conversion between logarithmic and fixed point, conversion circuits 
were realized as look-up tables and synthesized using the same process as FIR filter. Only 
conversion between LNS and fixed point is shown since FIR filters are typically implemented using 
fixed point. The results 
are shown in Table~\ref{tab:conv_ckt} and it can be seen that the cost of realizing FIR filters 
using LNS is still lower than using fixed point, albeit with reduced margins, including conversion 
results, up to word length of $Q(4,4)$. Optimized circuits for conversion have been proposed in the 
literature by Nam et al. in \cite{Nam2008J}, Zhu et al. in \cite{Zhu2016J} and others which can 
help in reducing this cost.

\section{Conclusion}
\label{sec:conc}
Logarithmic number system (LNSs) are a valuable alternative to
floating point for embedded systems and domain-specific hardware
accelerators. Existing LNS research has almost exclusively used
base-$2$. Although other bases, such as base-$\sqrt{2}$ or
base-$\sqrt[4]{2}$ have been explored, we show that these are simple
aliases of base-$2$, which can be constructed by moving the binary
point of the fixed-point exponent in a base-$2$ LNS.

In this paper we show that non-base-$2$ LNSs can have significant
advantages, particularly for short word lengths from $5-16$ bits.  We
first define a numeric error for LNS that measure \LULP~as an absolute
error in the fixed-point exponent of the LNS, which corresponds to a
multiplicative error in the corresponding real domain value.
We show that errors in conversion to LNS can be reduced dramatically by
choosing an appropriate base and scaling factor, at least for the
distribution of FP numbers. In our experiments, we reduce the typical
conversion error from $0.20 - 0.25$ \LULP~to less than $0.02$ \LULP.

We note that LNSs with a fixed word size are not closed under addition
or subtraction, which leads to rounding errors in these operations.
We also show that for a given LNS word size, some LNS bases lead to
inherently lower average rounding errors. When rounding to nearest,
the maximum arithmetic error for non-special values is 0.5 \LULP, and
we expect the average to be around 0.25 \LULP. For example, in 5-bit
(i.e. Q(2,2)) LNS, we show that base-$2$ leads to average errors of
0.260 (add) and 0.254 (subtract) \LULP. In contrast, base-$1.417$ gives
an average error of just 0.227(add) and 0.172 (subtract) \LULP.
Over millions of operations, even small differences in rounding
errors can accumulate to large errors.

Although LNS add/subtract are often implemented with low-latency ROM
lookup tables, we show that for small word sizes implementing these
functions in dedicated logic is also efficient. The choice of base
affects the truth table of the functions computed for LNS
add/subtract, and can therefore affect the area and latency of the
add/subtract units. By appropriate choice of base, the logic
area of the \phiplus and \phiminus logic can be reduced by around
2$\times$.

We also analyzed mixed precision tables, where the values stored in
the tables were of higher precision than the input. This can reduce
the errors for low precision arithmetic. The building blocks of adder, subtractor and multiplier 
are implemented using logarithmic, fixed and floating point numbers and extended to a complete FIR 
filters and we showed that LNS can be beneficial from all aspects of area, delay and power for very 
short word lengths.

We conclude that the choice of base affects the conversion
errors, arithmetic errors, and hardware implementation of arithmetic
units on several dimensions and that LNS can be used to implement different systems, like the FIR 
filter. Different bases allow different trade-offs
between these goals. Although prior research in LNS has generally
assumed base-$2$, our results show that other bases can reduce
arithmetic and conversion errors while simultaneously reducing
the area and latency of add/subtract and FIR units. Thus, when designing
a low-precision LNS for custom hardware, one should consider
looking beyond base-$2$.

\section*{Acknowledgments}
	This material is based upon work supported, in part, by Science Foundation 
	Ireland under Grant No. 13/RC/2094\_P2 and, in part, by the 
	European Union's Horizon 2020 research and innovation programme under the Marie 
	Sk\l odowska-Curie grant agreement and Grant No. 754489. Any opinions, findings, 
	and conclusions or recommendations expressed in this material are those	of the author and do not 
	necessarily reflect the views of the Science Foundation Ireland and European Union's Horizon 2020 
	programme.
	
	We also extend our gratitude to Dr. Oscar Gustafsson of Link\"oping University, Sweden and 
	Institute of Technology Carlow, Carlow, Ireland for their support.

%\bibliographystyle{unsrt}
%%% Bibliography
%\bibliography{IEEEabrv,../../TCD_Conf_abrv,../../TCD_Journal_abrv,../../TCD_bib}

\end{document}